\documentclass{article}
\usepackage{epsfig,amsfonts}
\usepackage{amsmath}
\usepackage{amsthm}
\usepackage{graphicx,color}
\usepackage{pgf}
\usepackage{url}
\usepackage{lineno}
\usepackage{mathtools}
\usepackage{adjustbox}
\usepackage{algpseudocode}
\usepackage{setspace} % use this to set more space between lines in algorithm outlines

\newcommand{\reals}{\mathbb{R}}%{{\rm I\!R}}

\newtheorem{algorithm}{Algorithm}[section]

\usepackage{lineno}
\usepackage{mathtools}
\usepackage{adjustbox}
\usepackage{algpseudocode}
\usepackage{setspace} % use this to set more space between lines in algorithm outlines

\newcommand\Tstrut{\rule{0pt}{2.6ex}}       % Top strut
\newcommand{\resetcounters}{\setcounter{equation}{0} \setcounter{figure}{0}
 \setcounter{table}{0}}

\begin{document}

\title{Avoiding breakdown in incomplete factorizations in low precision arithmetic
}

\author{Jennifer Scott\thanks{
STFC Rutherford Appleton Laboratory,
Harwell Campus, Didcot, Oxfordshire, OX11 0QX, UK
and School of Mathematical, Physical and Computational Sciences,
University of Reading, Reading RG6 6AQ, UK.
Correspondence to: {\tt jennifer.scott@stfc.ac.uk.}}
\and Miroslav T\accent23uma\thanks{
Department of Numerical Mathematics, Faculty of Mathematics and Physics,
Charles University, Czech Republic. {\tt mirektuma@karlin.mff.cuni.cz.} 
}
}

\maketitle

\begin{abstract}

The emergence of low precision floating-point
arithmetic in computer hardware has led to a resurgence of interest in the use of mixed precision
numerical linear algebra. For linear systems of equations,
there has been renewed enthusiasm for mixed precision variants of 
iterative refinement. We consider the iterative solution of large sparse systems
using incomplete factorization preconditioners.
The focus is on the robust computation of such preconditioners 
in half precision arithmetic and employing them to solve
symmetric positive definite systems to higher precision
accuracy; however, the proposed ideas can be applied more generally.
Even for well-conditioned problems,  incomplete factorizations 
can break down when small entries occur on the diagonal during the factorization.
When using half precision arithmetic, overflows are an additional possible source of breakdown.
We examine how breakdowns can be avoided and 
we implement our strategies within new half precision Fortran sparse
incomplete Cholesky factorization software. Results are reported
for a range of problems from practical applications. 
These demonstrate that, even for highly ill-conditioned problems, half precision preconditioners can potentially
replace double precision preconditioners, although unsurprisingly this may be at the cost
of additional iterations of a Krylov solver.
\end{abstract}

\maketitle

\section{Introduction}

We are interested  in solving sparse linear systems $Ax = b$, where $A \in \reals^{n \times n}$ 
is nonsingular and $x, b \in \reals^{n}$.
The majority of algorithms for solving such systems fall into two main categories:
direct methods and iterative methods. Direct methods
transform $A$ using a finite sequence of elementary transformations
into a product of simpler sparse matrices
in such a way that solving linear systems of equations with the factor matrices is 
comparatively easy and relatively inexpensive. For example, for a general nonsymmetric matrix,
$A = PLUQ$, where $L$ is a lower triangular matrix, $U$ is an upper
triangular matrix, and $P$ and $Q$ are permutation matrices
chosen to preserve sparsity in the factors and ensure
the factorization is stable. Direct methods, when properly implemented, are robust, frequently very fast, and can be 
confidently used as black-box solvers for computing  solutions
with predictable accuracy.
%, which is typically {at the level of double precision}.
However, they require  significant expertise to implement efficiently
(particularly in parallel). They can also need large amounts of memory
(which increases nonlinearly with the size and density of $A$) and the matrix factors normally contain
many more nonzero entries than $A$; these extra entries
are termed the fill-in and much effort goes into trying to minimise the amount of fill-in.

By contrast, iterative methods compute a sequence of approximations
$ {x}^{(0)}, {x}^{(1)}, {x}^{(2)}, \ldots $
that (hopefully) converge to the solution in an acceptable number of iterations.
The number of iterations (and whether or not convergence occurs at all) 
depends on  ${x}^{(0)}$, $A$ and $b$ as well as the required accuracy
in $x$. Basic implementations of iterative solvers are relatively straightforward 
as they only use the sparse matrix $A$ indirectly,
through matrix-vector products and, most importantly, their memory demands are limited
to a (small) number of vectors of length $n$, making them attractive for
very large problems as well as for problems where $A$ is not available explicitly.
However,  preconditioning is usually essential to enhance convergence of the iterative method.
Preconditioning seeks to transform the system into one that is more tractable
and from which the required solution of the original system can easily be recovered.
Determining and computing effective preconditioners is highly problem
dependent and generally very challenging.
Algebraic preconditioners that are built using an incomplete factorization of $A$ in which entries
that would be part of a complete factorization are dropped are frequently used,
especially when the underlying physics of the problem is difficult to exploit. Such
preconditioners can be
employed within more sophisticated methods; for example, to precondition subdomain solves 
in domain decomposition schemes or as smoothers in multigrid methods.

The performance differences for computing and communicating in different precision formats has led to
a long history of efforts to enhance numerical algorithms by combining
precision formats. The goals of mixed-precision
algorithms include accelerating the computational time by including the use of
lower-precision formats while obtaining high-precision accuracy
of the output and, by reducing the memory requirements, extending
the size of problems that can be solved as well as potentially lowering energy usage.
Numerical linear algebra software, and linear system solvers 
in particular, traditionally use double precision (64-bit) arithmetic, although some
packages (including the  BLAS and  LAPACK routines and some sparse solvers, such as those in the HSL
mathematical software library \cite{hsl:2023}) have always also offered single precision  (32-bit) versions.
In the mid 2000s, the speed difference between single and double precision arithmetic on
what were then state-of-the-art architectures, notably
Sony/Toshiba/IBM Cell processors (see, for example, \cite{Buttari:2007,kudo:2007,lllk:2006}), led to studies into the
feasibility of  factorizing a matrix in single precision
and then using the factors as a preconditioner for a simple iterative method
to regain higher precision accuracy \cite{ardu:2009,Buttari:2008,Buttari:2007a}. 
Hogg and Scott~\cite{hosc:2010} extended this to 
develop a mixed precision sparse symmetric indefinite solver.
More recently, the potential for employing single precision arithmetic in solving sparse systems to double
precision accuracy using multiple cores has been considered by Zounon et al~\cite{zhlt:2022}.
For such systems, the integer data used to hold the sparse structures of the matrix and its factors 
is independent of the precision and memory savings come only
from the real factor data. This offers the potential to increase the size of problem
that can be tackled and, in the incomplete factorization case,
to allow more entries to be retained in the factors, which can
result in a higher quality preconditioner, leading to savings in the total solution time.

In the past few years, the emergence of lower precision arithmetic in hardware has
led to further interest in mixed precision algorithms. A key difference
compared to  earlier work is the use of half precision (16-bit)
arithmetic. A comprehensive state-of-the-art survey of work on mixed precision numerical 
linear algebra routines (including an extensive bibliography)
is given in \cite{hima:2022} (see also \cite{abdel:2021}).
In particular, there have been important ideas and theory
on mixed precision iterative refinement methods that employ
the matrix factors computed in low precision as a preconditioner
to recover higher precision accuracy \cite{ames:2021,cahi:2017,cahi:2018}.
The number of entries in the factors of a sparse matrix $A$ is typically much greater than
in $A$ and so, unlike in the dense case, the overhead of keeping a high precision copy
of $A$ for computing the residual in the refinement process is small.
Amestoy et al~\cite{ames:2023} investigate the potential of mixed precision 
iterative refinement to enhance methods for sparse systems based on 
a particular class of approximate sparse factorizations. 
They employ the well-known parallel sparse direct solver MUMPS~\cite{amdd:2000}, which is able
to exploit block low-rank  factorizations
and static pivoting to compute approximate factors. The reported results in \cite{ames:2023}
are restricted to combining single and double precision arithmetic because,
in common with all other currently available
sparse direct solvers, MUMPS does not support the use of half precision arithmetic. Developing 
an efficient half precision solver is not straightforward, requiring 16-bit 
versions of the dense linear algebra routines that provide the building
blocks behind direct solvers. Such routines are
becoming available (for example, dense matrix-matrix multiplication 
in half precision is now supported by NVIDIA’s cuBLAS library), making the future
development of efficient software potentially more feasible.

Higham and Pranesh~\cite{hipr:2021} focus on symmetric positive definite linear systems.
They compute a Cholesky factorization using low precision arithmetic and employ the factors as preconditioners 
in GMRES-based and CG-based iterative refinement.
While they are interested in the sparse case, their MATLAB experiments
(which simulate low precision using their {\tt chop} function \cite{hipr:2019})
store the sparse test examples as dense matrices and their Cholesky
factorizations are computed using dense routines. The
reported theoretical and numerical results demonstrate the potential for
low precision (complete) factors to be used to obtain high precision accuracy.
Most recently, Carson and Khan~\cite{cakh:2023} have considered
using sparse approximate inverse preconditioners (SPAI) that are
based on Frobenius norm minimization \cite{grhu:97}.
They again propose computing the preconditioner using low precision and then 
employing GMRES-based iterative refinement and report
MATLAB experiments using the {\tt chop} function.
Mixed precision has also been investigated for
multigrid methods (see, for example,  the error analysis of McCormick et al~\cite{mcbt:2021}, 
who observed that different levels
in the grid hierarchy should use different precisions). 

Our emphasis is on low precision incomplete factorization preconditioners,
combined with Krylov subspace-based iterative refinement.
Although our ideas can be used for general sparse linear systems,
we focus on the sparse symmetric positive definite
case. We use half precision arithmetic to construct incomplete Cholesky factorizations that are then
employed as preconditioners to recover double precision accuracy.
Our primary objective is to show that for a range of problems (some of which
are highly ill-conditioned) it is possible to successfully obtain and use low precision incomplete factors.
We consider the potential sources of overflow during the incomplete factorization
and look at how to safely detect and prevent overflow. 
{We follow a number of others working on the development
of numerical linear algebra algorithms in mixed precision in
performing experiments that aim  to explore the 
feasibility of the ideas by simulating half precision (see, for example,  \cite{hipr:2019,cahp:2020a,cakh:2023}).
We want the option to experiment with sparse problems
that may be too large for MATLAB and have 
chosen to develop our software in Fortran.}
Half precision and double precision
versions are tested on systems coming from practical applications.

This paper makes the following contributions:
\begin{description}
\item{(a)} it considers the practicalities of computing classical incomplete
factorizations using low precision arithmetic, in particular, the 
safe prediction of possible overflows;
\item{(b)} it successfully employs global modifications
combined with a simple prescaling to prevent breakdowns during the factorization;
\item{(c)} it develops level-based incomplete Cholesky factorization software
in half precision (written in Fortran);
\item{(d)} it demonstrates the potential for robust half precision incomplete factorization preconditioners
to be used to solve ill-conditioned symmetric positive definite systems.
\end{description}

The rest of the paper is organised as follows. In Section~\ref{sec:incomplete}, we briefly recall
incomplete factorizations of sparse matrices and consider the challenges that incomplete
Cholesky factorizations can face, particularly when low precision arithmetic is employed. In Section~\ref{sec:refine},
we summarise basic mixed precision iterative refinement algorithms before presenting
numerical results for a range of problems coming
from practical applications in Section~\ref{sec:experiments}. Concluding remarks
and possible future directions are given in Section~\ref{sec:conclusions}.

\medskip\noindent
{\bf Terminology.} 
We use the term high precision for precision
formats that provide high accuracy at the cost of a larger
memory volume (in terms of bits) and low precision to refer
to precision formats that compose of fewer bits (smaller
memory volume) and provide low(er) accuracy. Unless
stated otherwise, we mean IEEE double precision (64-bit) when
using the term high precision (denoted by fp64) and  the 1985 IEEE standard 754 half precision (16-bit)
when using the term low precision (denoted by fp16).  bfloat16 is another form of half-precision arithmetic
that was introduced by
Google in its tensor processing units and  formalized by Intel; we do not use it in this paper. 
Table~\ref{table:precisions} summarises  the parameters for different precision arithmetic. 
We  use $u_{16}$ , $u_{32}$ , $u_{64}$ to denote the unit roundoffs
in fp16, fp32, and fp64 arithmetic, respectively.
\begin{table}[htbp]
\begin{center}
\caption{Parameters for bfloat16, fp16, fp32, and fp64 arithmetic: the number of bits in the significand and
exponent, unit roundoff $u$, smallest positive (subnormal) number $x^s_{min}$ , smallest normalized positive
number $x_{min}$, and largest finite number $x_{max}$, all given to three significant figures.
${\dagger}$ In Intel's bfloat16 specification, subnormal numbers are not supported.
}
 \label{table:precisions}
\vspace{3mm}
\footnotesize
\begin{tabular}{cllllll}
\hline\Tstrut
&
Signif. &
 Exp. &
$ \;\;u $ &
 $\;x^s_{min}$  &
 $\;x_{min}$  &
 $\;x_{max}$  \\
 \hline\Tstrut
bfloat16  & 8 & 8 & $3.91 \times  10^{-3}$ & $\;\;\; \dagger$ & $1.18 \times  10^{-38}$ & $3.39 \times 10^{38}$ \\
fp16  & 11 & 5 & $4.88 \times 10^{-4}$ & $5.96 \times 10^{-8}$ & $6.10 \times 10^{-5}$ & $6.55 \times 10^4$ \\
fp32 & 24 & 8 & $5.96 \times 10^{-8}$ & $1.40 \times 10^{-45}$ & $1.18 \times 10^{-38}$ & $3.40 \times 10^{38}$ \\
fp64 & 53 & 11 & $1.11 \times 10^{-16}$ &$4.94 \times 10^{-324}$ & $2.22 \times 10^{-308}$ & $1.80 \times 10^{308}$ \\

\hline
\end{tabular}
\end{center}
\end{table}

\section{Incomplete factorizations}\label{sec:incomplete}
In this section, we briefly recall incomplete factorizations of
sparse matrices and then discuss how breakdown can happen, particularly
when using low precision.

\subsection{A brief introduction to incomplete factorizations}

For a general sparse matrix $A$, the incomplete factorizations that we are interested in are of the form $A \approx LU$, where
$L$ and $U$ are sparse lower and upper triangular matrices, respectively (for simplicity of notation,
the permutations $P$ and $Q$ that are for preserving sparsity
in the factors are omitted). The computed incomplete factors are used to define
the preconditioner; the preconditioned linear system is $U^{-1} L^{-1}A = U^{-1} L^{-1}b$.
If $A$ is a symmetric positive definite (SPD) matrix then $A = LL^T$.
There are three main classes of incomplete factorization preconditioners.
Firstly, threshold-based $ILU(\tau)$ methods in which the locations of permissible fill-in
in the factors are determined in conjunction with the numerical factorization of $A$; entries of
the computed factors that are smaller than a prescribed threshold are dropped e.g., \cite{munk:80}.
Secondly, memory-based $ILU(m)$ methods in which the amount of memory available for
the incomplete factorization is prescribed and only the largest entries are retained
at each stage of the factorization e.g., \cite{jopl:95}.
Thirdly, structure-based $ILU(\ell)$ methods in which an initial symbolic phase
uses a hierarchy of sparsity structures and
determines the location of permissible entries using only the sparsity pattern of $A$ e.g., \cite{watt:81}.
The memory requirements for the incomplete factors are then determined before
the numerical factorization is performed.
The simplest such approach (for which the symbolic phase is trivial) is an $ILU(0)$ factorization (or $IC(0)$ in the 
SPD case) that limits entries in the incomplete factors
to positions corresponding to entries in $A$ (no fill-in is permitted). $ILU(0)$ preconditioners are frequently
used for comparison purposes when assessing the performance of other approaches.

The different approaches have been developed, modified and refined over many years. 
Variants have been proposed that combine
the ideas and/or employ them in conjunction with discarding
entries in $A$ (sparsification) before the factorization commences.
{For an introduction to ILU-based preconditioners
and details of possible variants, we recommend 
\cite{chvd:97} and Chapter 10 of the book \cite{saad:03}, 
while \cite{sctu:2011} provides a brief historical overview and further
references.

Algorithm~\ref{alg:ic_generic} outlines a basic (right-looking) incomplete
Cholesky (IC) factorization of a SPD matrix $A= \{a_{ij}\}$. 
Here all computations are performed in the same precision as $A$.
The algorithm assumes a target sparsity pattern ${\mathcal S}\{L\}$ for $L=\{l_{ij}\}$ is provided, where
$${\mathcal S}\{L\}  = \{ (i,j) \, | \ l_{ij} \neq 0, \, 1 \le  j \le i \le n \}.$$
{It is assumed that ${\mathcal S}\{L\}$ includes the diagonal entries.}
Modifications can be made to incorporate threshold dropping strategies
and to determine ${\mathcal S}\{L\}$ as the method proceeds.
At each major step $k$,  outer product updates are applied to the part of the matrix
that has not yet been factored (Lines 7--11).

\medskip
\hrule\vskip-2mm
\begin{algorithm}{\bf Basic right-looking IC factorization}\\
\textbf{Input:}  SPD matrix $A$ and a target sparsity pattern ${\mathcal S}\{L\}$  \\
\textbf{Output:} Incomplete Cholesky factorization  $A \approx  L L^T$.
\label{alg:ic_generic}
\smallskip\hrule\smallskip
\setstretch{1.17}\begin{algorithmic}[1]
\State Initialize $l_{ij} = a_{ij}$ for all $(i,j) \in {\mathcal S}\{L\}$
\For{$k=1:n$}\Comment{Start of $k$-th major step}
  \State $ l_{kk}\leftarrow (l_{kk})^{1/2}$\Comment{Diagonal entry is the pivot}
  \For{$i=k+1:n$ such that $(i,j) \in {\mathcal S}\{L\}$}
    \State    $l_{ik} \leftarrow l_{ik} / l_{kk}$ \Comment{Scale pivot column $k$ of the incomplete factor by the pivot}
  \EndFor \Comment{Column $k$ of $L$ has been computed}
  \For{$j=k+1:n$  such that $(j,k) \in {\mathcal S}\{L\}$}
    \For{$i=j:n$ such that $(i,j) \in {\mathcal S}\{L\}$ }
     \State    $l_{ij} \leftarrow l_{ij} - l_{ik} l_{jk}$ \Comment{Update operation on column $j$}
   \EndFor
  \EndFor
\EndFor
\end{algorithmic}
\end{algorithm}
\null\vskip-6mm\hrule\smallskip
\medskip

\subsection{Breakdown during incomplete factorizations}
\label{sec:breakdowns}
For arbitrary choices of the sparsity pattern ${\mathcal S}\{L\}$,
the incomplete Cholesky factorization exists if $A$ is an M-matrix or a H-matrix
with positive diagonal entries \cite{meva:77,mant:80}. But for a general SPD
matrix, there is no such guarantee and an incomplete factorization algorithm can 
(and frequently does) break down. There are three places where breakdown
can occur. We refer to these as problems B1, B2, and B3 as follows.
\begin{itemize}
\item B1: The computed diagonal entry $l_{kk}$ (which is termed the pivot
at step $k$) may be unacceptably small or negative. 
\item B2: The scaling $l_{ik} \leftarrow l_{ik} / l_{kk}$ may overflow.
\item B3: The update $l_{ij} \leftarrow l_{ij} - l_{ik} l_{jk}$ may overflow.
\end{itemize}

It is crucial for a robust implementation to detect the possibility of overflow
before it occurs (otherwise, the code will simply crash). Tests 
for potential breakdown must themselves only use operations that cannot overflow.
We say that an operation is  {\it safe} in the precision being used
if it cannot overflow; otherwise it is  {\it unsafe}. 

Safe detection of problem B1 is straightforward as
we  simply need to check at Line 3 of Algorithm~\ref{alg:ic_generic}
that the pivot satisfies $l_{kk} \ge \tau > x_{min} $, where the  threshold parameter $\tau$
depends on the precision used. Based on our experience with practical problems, in our reported
results in Section~\ref{sec:experiments}, we use $\tau = 10^{-5}$ for 
half precision factorizations and,
for double precision factorizations,  $\tau = 10^{-20}$ (note that preconditioner quality is
not sensitive to the precise choice of $\tau$). 

Problem B2 can happen during the column scaling at Line 5 of Algorithm~\ref{alg:ic_generic}.
Let $a$ be the entry of largest absolute
value in the current pivot column $k$ and 
let $d$ denote the current pivot $l_{kk}$ (set at Line 3).
If $a \le x_{max}$ and $1 \le d \le x_{max}$ then
$d \ge a/x_{max}$ and it is safe to compute $a/d$ (and thus safe to scale column $k$).

Problem B3 can occur during the update operations at Line 9 of Algorithm~\ref{alg:ic_generic}.
Algorithm~\ref{alg:safe_update} can be used to perform safe updates. Here the scalars and the computation
are in the working precision.
\medskip
\hrule\vskip-2mm
\begin{algorithm}{\bf Safe test for a safe update operation (detecting B3 breakdown)}\\
\textbf{Input:}  Scalars $a,b,c$ such that $|a|,|b|,|c| \le x_{max}$.  \\
\textbf{Output:} Either $v=a-bc$ or $flag = -3$ (unsafe to perform update).
\smallskip\hrule\smallskip
\setstretch{1.17}\begin{algorithmic}[1]
    \State  $flag = 0$
    \If {$|b| \le 1$ or  $|c| \le 1$} \Comment{$x_{max}/|c|$ cannot be tested if $|c| < 1$}
      \State    $w = bc$
    \ElsIf {$|b| \ \le \ x_{max}/|c|$} 
      \State    $w = bc$
    \Else
       \State $flag = -3$ and {\bf return} \Comment{Unsafe to compute $bc$}
    \EndIf
    \If  {$a \ge 0$ }  
      \If  {$w \ge 0$   or ($x_{max}-a  \ge  -w$)}
         \State $v = a - w$
      \Else
         \State $flag = -3$ and {\bf return}  \Comment{Unsafe to perform subtraction}
      \EndIf
     \Else
       \If  {$w < 0$   or ($x_{max}+a  \ge  w$)}
          \State $v = a - w$
       \Else
         \State $flag = -3$ and {\bf return}  \Comment{Unsafe to perform subtraction}
       \EndIf
     \EndIf
\end{algorithmic}
\label{alg:safe_update}
\end{algorithm}
\null\vskip-6mm\hrule\smallskip
\medskip

Algorithm~\ref{alg:ic_generic2} presents the basic IC factorization algorithm with the inclusion of safe
checks for possible breakdown. Here the matrix $A$ and the computation are in the working precision.
In practice, when using single and double precision  the tests in Lines 8--11 and 16--20 are omitted
but B1 can occur when using any arithmetic. Note that the safe test 
for B3 breakdown does not have to be applied to individual scalar entries
but can be applied to each column $j \ge k+1$.
\newpage%\medskip
\hrule\vskip-2mm
\begin{algorithm}{\bf Basic IC factorization with safe checks for breakdown} \\
\textbf{Input:}  SPD matrix $A$, a target sparsity pattern ${\mathcal S}\{L\}$ and parameter $\tau>0$ for
testing small entries. \\
\textbf{Output:} Either $flag < 0$ (breakdown occurred) or $A \approx  L L^T$ with $L$ lower triangular
\smallskip\hrule\smallskip
\label{alg:ic_generic2}
\setstretch{1.17}
\begin{algorithmic}[1]
\State Initialize $l_{ij} = a_{ij}$ for all $(i,j) \in {\mathcal S}\{L\}$ 
\State Set $flag = 0$
\For{$k=1:n$}\Comment{Start of $k$-th major step}
\If{$l_{kk} < \tau$}
   \State Set $flag = -1$ and  {\bf return}\Comment{B1 breakdown}
  \EndIf
    \State $ l_{kk}\leftarrow (l_{kk})^{1/2}$
    \State $a= \max_{i=k+1:n} \{ |l_{ik}|: (i,k) \in {\mathcal S}\{L\}\}$
    \If{ $l_{kk} \ge 1$ or $l_{kk} \ge a/x_{max}$} \Comment{Note that if $l_{kk} \ge 1$ then $a$ does not need to be computed}
  \For{$i=k+1:n$ such that $(i,k) \in {\mathcal S}\{L\}$}
    \State  $l_{ik} \leftarrow l_{ik} / l_{kk}$ \Comment{Perform safe scaling}
  \EndFor\Comment{Column $k$ of $L$ has been computed}
  \Else
         \State Set $flag = -2$ and  {\bf return}\Comment{B2 breakdown}
  \EndIf
  \For{$j=k+1:n$ such that $(j,k) \in {\mathcal S}\{L\}$}
%      \If{$ \max_{i=k+1:n} \{l_{ik}^2: (i,k) \in {\mathcal S}\{L\}\} + \max_{i=j:n} \{|l_{ij}|: (i,j) \in {\mathcal S}\{L\}\} > x_{max}$}
%      \State Set $a= \max_{i=j:n} \{|l_{ij}|: (i,j) \in {\mathcal S}\{L\}\}$
%      and $b =c =\max_{i=k+1:n} \{|l_{ik}|: (i,k) \in {\mathcal S}\{L\}\}$
      \State Use Algorithm~\ref{alg:safe_update} to test for safe update %\cblue{Only need to test $b \ \le \ x_{max}/b$???}
       \If{ $flag = -3$}  
       \State{\it \bf return}\Comment{B3 breakdown}
    \EndIf
    \For{$i=j:n$ such that $(i,j) \in {\mathcal S}\{L\}$ }
     \State    $l_{ij} \leftarrow l_{ij} - l_{ik} l_{jk}$\Comment{Perform safe update operation}
   \EndFor
  \EndFor\Comment{Column $j$ of $L$ has been updated}
\EndFor
\end{algorithmic}
\end{algorithm}
\null\vskip-6mm\hrule\smallskip
\medskip

Observe that the occurrence of underflows 
when using fp16 arithmetic does not prevent the computation of the incomplete factors,
although underflows could lead to a loss of information
that  affects the preconditioner quality.
However, provided the problem has been well scaled, the dropping strategy 
used within the incomplete factorization
normally has more influence on the computed factors than underflows do.
A subnormal floating-point number is a nonzero number with absolute value less than that
of the smallest normalized number. Floating-point operations on subnormals can be very slow, because they often require
extra clock cycles, which introduces a high overhead.
If an off-diagonal factor entry is subnormal, it can again be replaced by zero without significantly
affecting the preconditioner quality.

\subsection{Prescaling of the system matrix}

Sparse direct solvers typically prescale $A$ and factorize $\widehat A = S_1^{-1}AS_2^{-1}$, where $S_1$ and $S_2$
are diagonal scaling matrices ($S_1= S_2 = S$ in the symmetric case). Because no single choice of scaling
always results in the best performance
(in terms of the time, factor sizes, memory requirements and data movement), several
possibilities (with different associated costs) are typically offered
to allow a user to experiment and select the optimum for their application.
For incomplete factorizations, prescaling can reduce the incidence of breakdowns. 
This is illustrated for SPD matrices
in \cite{sctu:2014a}, where experiments (in double precision arithmetic) show that 
the cheap  scaling in which the entries in column $j$ of $A$ are normalised
by the 2-norm of column $j$  (so that the absolute values
 of the entries of the scaled matrix are all less than 1) is generally a good choice. 
However, scaling alone cannot guarantee to prevent breakdown. If breakdown does
happen then modifications need to be made to the scaled matrix that is being factorized,
either before or during the factorization; this
is discussed in the next subsection.

In their work on matrix factorizations in low precision, Higham, Pranesh and Zounon~\cite{hipz:2019} prescale
$A$ using the working precision and then round to fp16 once all
the absolute values of all the entries of the matrix are at most $x_{max}$ 
(see also \cite{hipr:2021,hima:2022}).
They employ an equilibration approach \cite{knru:2014} and then multiply
each entry of the scaled matrix by $\mu = \theta x_{max}/\gamma$, where $\gamma$ is the 
entry of largest absolute value in the scaled matrix $\widehat A$ and $\theta \in (0,1)$ is a parameter chosen with
the objective of limiting 
the possibility of overflow during the factorization. 
Our experiments have found, in the incomplete factorization case, that it is sufficient to use
an inexpensive $l_2$-norm scaling  and we set $\mu = 1$.

\subsection{Global modifications to prevent breakdown}
\label{sec:global mods}
Local diagonal modifications were first described in the 1970s by  
Kershaw~\cite{kers:78}.  The idea is to simply modify an individual diagonal entry of
$A$ during the factorization if it is found to be too small (or negative) to use as a pivot, that is,
at the Lines 4--6 of Algorithm~\ref{alg:ic_generic2}, instead of terminating the factorization, the  pivot is
perturbed by some positive quantity so that it is at least $\tau$. 
While such local modifications are inexpensive to implement within a right-looking factorization algorithm, 
it is frequently the case that even a small number 
of modifications can result in a
poor preconditioner. Consider, for example, the extreme case that,
at the start of the major step $k$ of Algorithm~\ref{alg:ic_generic2} the diagonal entry 
$l_{kk}$ is equal to zero and the absolute values of all the remaining
entries of column $k$ are $x_{max}$. A local modification that replaces $l_{kk}$ by
some chosen value less than $1$ prevents B1, but the corresponding column then does not scale 
(entries overflow) so that B1 is transferred into a B2 problem.

Global strategies are generally more successful
in terms of the quality of the resulting preconditioner (see, for example, \cite{benz:02}
and the theoretical  and numerical results given in \cite{limo:99,mant:80,sctu:2014a}).
When an incomplete factorization breaks down, a straightforward strategy is to select a shift $\alpha>0$,
replace the scaled matrix $\widehat A$ by $\widehat A + \alpha I$ 
($I$ is the identity matrix) and restart the factorization. 
This is outlined in Algorithm~\ref{alg:trial_error}. The factors of the 
shifted matrix are used to precondition $\widehat A$. 
If $A_D$ and $A_E$ are, respectively, the diagonal and off-diagonal parts of $\widehat A$,
then there is always some  $\alpha$ for which  $(1 +\alpha) A_D + A_E$
is diagonally dominant. Provided the target sparsity pattern
of the incomplete factors contains the positions of the 
diagonal entries, then it can be shown that the incomplete factorization of this shifted matrix 
does not break down \cite{mant:80}. Diagonal dominance is sufficient
for avoiding breakdown but it is not a necessary condition and an
incomplete factorization may be breakdown free for much smaller values of $\alpha$. An appropriate  $\alpha$
is not usually known a priori: too large a value may harm the quality of the incomplete factors
as a preconditioner  and too small a value will not prevent breakdown,
necessitating more than one restart, with a successively larger $\alpha$. Typically, the shift is doubled after
a breakdown, although more sophisticated strategies are sometimes used (see, for example, \cite{sctu:2014a}).

\medskip
\hrule\vskip-2mm
\begin{algorithm}{\bf Shifted incomplete IC factorization} \\
\textbf{Input:}  SPD matrix $A$, scaling matrix $S$, and initial shift $\alpha_S \ge 0$, all in precision $u$ \\
\textbf{Output:} Shift $\alpha \ge 0$ and incomplete 
Cholesky factorization $S^{-1}AS^{-1}+ \alpha I \approx LL^T$ in precision $u_l \ge u$
\setstretch{1.17}
\begin{algorithmic}[1]
\State $\alpha_0 = 0$
\State $\widehat A = S^{-1}AS^{-1} $ in precision $u$
\State if $u_l \ne u$ then $\widehat A^{(\ell)} = fl(\widehat A)$ \Comment{Squeeze the scaled matrix into precision $u_l$}
\For{$i = 0,1,2, \ldots$}
 \State $\widehat A^{(\ell)} +\alpha_i I\approx L L^T$ in precision $u_l$\Comment{Use Algorithm~\ref{alg:ic_generic2}}
 \State if successful then set $\alpha = \alpha_i$ and {\bf return}\Comment{This is the case of no breakdown}
 \State $\alpha_{i+1} =   \max(2 \alpha_i, \ \alpha_S)$\Comment{Breakdown detected so increase the shift and restart}
 \EndFor
\end{algorithmic}
\label{alg:trial_error}
\end{algorithm}
\null\vskip-6mm\hrule\smallskip
\medskip

When using fp16 arithmetic for the factorization ($u_{\ell} = u_{16}$),  Line 3
in Algorithm~\ref{alg:trial_error}  ``squeezes'' the scaled matrix
$\widehat A$ from the working precision $u$ into half precision. The squeezed
matrix $\widehat A^{(\ell)}$ is factorized using Algorithm~\ref{alg:ic_generic2} with 
everything in precision $u_l$.
Our experiments on SPD matrices confirm that, provided we prescale $A$,
the number of  times we must increase $\alpha$ and restart
is generally small but it depends
on the IC algorithm and the precision (see the statistics $nmod$ and $nofl$ is the tables of results in Section~\ref{sec:results}).

\subsection{Using the low precision factors}

Each application of an incomplete LU factorization preconditioner 
is equivalent to solving a system $LU v = w$ (or $LL^Tv = w$ in the SPD case). This involves
a solve with the lower triangular $L$ factor followed by a solve with the 
upper triangular $U$ (or $L^T$) factor; these are
termed forward and back substitutions, respectively.
Algorithm~\ref{alg:L solve} outlines a simple lower triangular solve.
\smallskip
\hrule\vskip-2mm
\begin{algorithm}{\bf Forward substitution: lower triangular solve $Ly= w$}\\
\textbf{Input:}  Lower triangular matrix $L$ with nonzero diagonal entries
and  right-hand side $w$. \\
\textbf{Output:} The solution vector $y$.
\smallskip\hrule\smallskip
\setstretch{1.15}\begin{algorithmic}[1]
 \State Initialise $y_j = w_j$, $1 \le j \le n$
\For{ $j=1:n$}
\If {$\;y_j \ne 0$} 
\State $y_j \leftarrow y_j /l_{jj}$
 \For{$i = j+1:n$}
   \If {$\;l_{ij} \neq 0$}
 \State $y_i \leftarrow y_i - l_{ij} y_j$
 \EndIf
 \EndFor
  \EndIf
 \EndFor
\end{algorithmic}
\label{alg:L solve}
\end{algorithm}
\null\vskip-6mm\hrule\smallskip
\medskip
If the factors are computed and stored in fp32 or fp64 arithmetic then
overflows are unlikely to occur in Algorithm~\ref{alg:L solve}. However,
if  half precision arithmetic $u_{\ell}$ is used for the factors and
the forward and back substitutions applied in precision 
$u_{\ell}$ then  overflows are much more likely at Lines 3 and 6. We can try and avoid this using 
simple scaling of the right hand side so that we solve $LU v = w/\|w\|_{\infty}$ (or $LL^T v = w/\|w\|_{\infty}$)
and then set $y = v \times \|w\|_{\infty}$ (see \cite{cahi:2018}).
Nevertheless,  as in problems B2 and  B3 above, overflows can still happen.
The safe tests for monitoring  potential overflow can
again be used. If detected,
higher precision arithmetic can be used for the triangular solves.

\section{LU and Cholesky factorization based iterative refinement}\label{sec:refine}
\label{sec:IR}

\subsection{LU-IR and Krylov-IR using low precision factors}

As already observed, 
a number of studies in the late 2000s looked at computing the (complete) matrix factors in single precision and then 
employing them as a preconditioner
within an iterative solver to regain double precision accuracy (e.g. \cite{ardu:2009,Buttari:2008,Buttari:2007a,hosc:2010}).
The simplest method is iterative refinement, which seeks to improve the
accuracy of a computed solution $x$ by iteratively repeating the following steps
until the required accuracy is achieved, the refinement stagnates, or a prescribed limit on the number
of iterations is reached.
\begin{enumerate}
\item Compute the residual $r = b - A  x$.
\item Solve the correction equation $Ad = r$.
\item Update the computed solution $  x \leftarrow  x + d$.
\end{enumerate}
A number of variants exist. For a general matrix $A$, the most common is
LU-IR, which computes the LU factors of $A$ in precision $u_{\ell}$, and then solves the correction equation
by forward and back substitution using the computed LU factors in precision $u_{\ell}$. 
The computation of the residual is performed in precision $u_r$ 
and the update is performed in the working precision $u$, with $u_r \le u \le u_{\ell}$. 
This is outlined in Algorithm~\ref{alg:lu-ir}.
Here and in Algorithm~\ref{alg:gmres-ir},  $A$ and  $b$ are held in the 
working precision $u$ and the computed solution is also in the working precision.
A discussion of stopping tests {(and many references) may be found in the book \cite{high:2002}
(see also the later papers \cite{bbdk:2009,dhkl:2006})}.

\medskip
\hrule\vskip-2mm
\begin{algorithm}{\bf LU-based iterative refinement using three precisions (LU-IR)}\\
\textbf{Input:}  Non singular matrix $A$ and vector $b$ in precision $u$, 
three precisions satisfying  $u_r \le u \le u_{\ell}$ \\
\textbf{Output:}  Computed solution of the system $Ax = b$ in precision $u$
\label{alg:lu-ir}
\smallskip\hrule\smallskip
\setstretch{1.17}\begin{algorithmic}[1]
\State  Compute the factorization $A = LU$ in precision $u_{\ell}$
\State  Initialize $x_1$ (e.g., by solving $LUx_1 = b$ using substitution in precision $u_{\ell}$)
\For{ i = 1 : itmax or until converged}\Comment{$itmax$ is the maximum iteration count}
\State  Compute $r_i = b - Ax_i$ in precision $u_r$; store $r_i$ in precision $u$
\State  Use the computed factors to solve $Ad_i = r_i$ by substitution in precision $u_{\ell}$; store $d_i$ in precision $u$
\State  Compute $x_{i+1} \leftarrow x_i + d_i$ in precision $u$
\EndFor
\end{algorithmic}
\end{algorithm}
\null\vskip-6mm\hrule\smallskip
\medskip

\medskip
\hrule\vskip-2mm
\begin{algorithm}{\bf Krylov-based iterative refinement using precisions (Krylov-IR)}\\
\textbf{Input:}  Non singular matrix $A$ and vector $b$ in precision $u$, a Krylov subspace method, and five precisions $u_r$,
$u_g$, $u_p$, $u$, $u_{\ell}$\\
\textbf{Output:}  Computed solution of the system $Ax = b$ in precision $u$
\label{alg:gmres-ir}
\smallskip\hrule\smallskip
\setstretch{1.17}\begin{algorithmic}[1]
\State  Compute the factorization $A = LU$ in precision $u_{\ell}$
\State  Initialize $x_1$ (e.g., by solving $LUx_1 = b$ using substitution in precision $u_{\ell}$)
\For{ i = 1 : itmax or until converged}\Comment{$itmax$ is the maximum iteration count}
\State  Compute $r_i = b - Ax_i$ in precision $u_r$; store $r_i$ in precision $u$
\State  Solve $U^{-1} L^{-1} Ad_i = U^{-1} L^{-1} r_i$ using the Krylov method in precision $u_g$, with
 $U^{-1} L^{-1} A $ performed in precision $u_p$; store $d_i$ in precision $u$
\State  Compute $x_{i+1} \leftarrow x_i + d_i$ in precision $u$
\EndFor
\end{algorithmic}
\end{algorithm}
\null\vskip-6mm\hrule\smallskip
\medskip

LU-IR can stagnate. In particular, if the LU
factorization is performed in fp16 arithmetic, then LU-IR is only guaranteed to reduce the solution
error if the condition number $\kappa(A)$ satisfies $\kappa(A) \ll 2 \times 10^3$.
To extend the range of problems that can be tackled,
Carson and Higham~\cite{cahi:2017} propose a variant that uses GMRES preconditioned
by the LU factors to solve the correction equation. This is outlined in Algorithm~\ref{alg:gmres-ir},
with the Krylov subspace method set to GMRES.
Carson and Higham use two precisions
$u= u_{\ell}$ and $u_r = u_g = u_p = u^2$; this was later extended to 
allow up to five precisions \cite{ames:2021,cahi:2018}. If the LU
factorization is performed in fp16 arithmetic 
and $u_g=u_p=u_{64}$, then  the solution
error is reduced by GMRES-IR provided $\kappa(A) \ll 3 \times 10^7$.
Note that Algorithm~\ref{alg:gmres-ir} requires two convergence tests and  stopping criteria; firstly, for the
Krylov method on Line 5 (inner iteration) and secondly, for testing the updated solution (outer iteration).

In the SPD case,  a natural choice is to choose the Krylov method to be the
conjugate gradient (CG) method. However, the supporting rounding
error analysis applies only to GMRES, because it relies on the backward stability of 
GMRES and preconditioned CG is not guaranteed to be backward stable \cite{gree:97a}. This
is also the case for MINRES. Nevertheless, the  MATLAB results presented in \cite{hipr:2021} suggest 
that in practice CG-IR generally works as well as GMRES-IR.

We observe that Arioli and Duff~\cite{ardu:2009} earlier proposed a 
simplified two precision variant of Krylov-IR
in which $itmax$ was set to 1.
They used single and double precision and employed restarted FGMRES~\cite{saad:93} as the Krylov solver.
This choice was based on their experience that FGMRES is more robust than GMRES \cite{adgs:2007}.
 Hogg and Scott~\cite{hosc:2010} subsequently developed a single-double precision 
solver for large-scale
symmetric indefinite linear systems; this Fortran code is available as {\tt HSL\_MA79}
within the HSL library \cite{hsl:2023}. It uses a single precision multifrontal
method to factorize $A$ (it computes a sparse LDLT factorization)
and then mixed precision iterative refinement (that is, LU-IR
with $u_r = u = u_{64}$ and $u_{\ell} = u_{32}$).
If iterative refinement stagnates, {\tt HSL\_MA79} employs restarted FGMRES to try and obtain double precision accuracy,
that is, a switch is automatically made within the code from LU-IR to Krylov-IR  with $x_1$ 
 in Line 2 of Algorithm~\ref{alg:gmres-ir} taken to be the current approximation
to the solution computed using LU-IR and $itmax = 1$ (see Algorithms 1--3 of \cite{hosc:2010}).

\subsection{Generalisation to low precision incomplete factors}

GMRES-IR can be modified by
replacing $U^{-1} L^{-1}$ with a preconditioner $M^{-1}$. 
Work on using  scalar Jacobi and ILU(0) preconditioners has been reported by  Lindquist et al. \cite{lild:2020,lild:2022},
with tests performed on a GPU-accelerated node combining single and double precision arithmetic. Loe et al. \cite{loe:2021}
present an experimental evaluation of multiprecision strategies for GMRES on GPUs
using block Jacobi and polynomial preconditioners. Carson and Khan~\cite{cakh:2023} use
a sparse approximate inverse preconditioner computed in low precisions and present MATLAB results.
In addition, Amestoy et al. \cite{ames:2023} use the option within the MUMPS 
solver to compute sparse factors in single precision using block low-rank  factorizations
and static pivoting and then employ them within GMRES-IR to recover double 
precision accuracy.

\medskip
\hrule\vskip-2mm
\begin{algorithm}{\bf Krylov-based iterative refinement with an incomplete factorization preconditioner using five precisions (IC-Krylov-IR)}\\
\textbf{Input:}  SPD matrix $A$ and vector $b$ in precision $u$, a Krylov subspace method, and five precisions $u_r$,
$u_g$, $u_p$, $u$ and $u_{\ell}$\\
\textbf{Output:} Computed solution of the system $Ax = b$ in precision $u$
\label{alg:gmres-ir-ic}
\smallskip\hrule\smallskip
\setstretch{1.17}\begin{algorithmic}[1]
\State  Compute an incomplete Cholesky factorization of $A$ in precision $u_{\ell}$
\Comment{Use Algorithm~\ref{alg:trial_error} to compute $L$}
\State  Initialize $x_1 = 0$
\For{ i = 1 : itmax or until converged}\Comment{$itmax$ is the maximum iteration count}
\State  Compute $r_i = b - Ax_i$ in precision $u_r$; store $r_i$ in precision $u$
\State  Solve $Ad_i = r_i$ using the preconditioned Krylov method in precision $u_g$, with
 preconditioning and products with $A$ performed in precision $u_p$; store $d_i$ in precision $u
 $\Comment{Computed factors used as the preconditioner}
\State  Compute $x_{i+1} \leftarrow x_i + d_i$ in precision $u$
\EndFor
\end{algorithmic}
\end{algorithm}
\null\vskip-6mm\hrule\smallskip
\medskip
Our interest is in using the IC factors of the SPD matrix $A$ 
as preconditioners within CG-IR and GMRES-IR. Algorithm~\ref{alg:gmres-ir-ic}
summarises the approach, which we call IC-Krylov-IR to emphasise
that an IC factorization is used (this is consistent with the
notation used in \cite{cakh:2023}). Note that if $itmax = 1$ then the algorithm
simply applies the preconditioned Krylov solver to try and achieve the requested accuracy.
Algorithm~\ref{alg:lu-ir} can be modified in a similar way to obtain what we will call
the IC-LU-IR method.

\section{Numerical experiments}\label{sec:experiments}

In this section, we investigate the potential effectiveness and reliability of half precision IC preconditioners.
Our test examples  are SPD matrices taken from the  SuiteSparse Collection; they are listed in Table~\ref{T:test problems}.
In the top part of the table are those we classify as being well-conditioned (those for which our
estimate $cond2$ of the 2-norm condition number is less than $10^7$) and 
in the lower part are ill-conditioned examples.
We have selected problems coming from a range of application areas and of
different sizes and densities. Many are initially poorly scaled and some
(including the first three problems in Table~\ref{T:test problems}) contain entries that overflow
in fp16 and thus prescaling of $A$ is essential.
We use the $l_2$ norm scaling (computed and applied in double precision arithmetic).
We have performed tests using
equilibration scaling (implemented using the HSL routine {\tt MC77}  \cite{Ruiz:2001,ruuc:2011})
and found that the resulting preconditioner
is of a similar quality; this is consistent with \cite{sctu:2014a}.
The right-hand side vector $b$ is constructed by setting
the solution $x$ to be the vector of 1's.
% Mirek is not using any preordering of A in the experiments (iorder = 0)

\label{sec:results}
\resetcounters
\begin{table}[htbp]
\caption{Statistics for our test examples. 
Those in the top half are considered to be well conditioned
and those in the lower half to be ill conditioned. 
$nnz(A)$ denotes the number of entries in the lower triangular part of $A$.
$normA$ and $normb$ are the infinity norms of $A$ and $b$. $cond2$
is a computed estimate of the condition number of $A$.
}
\label{T:test problems}\vspace{3mm}
{
\footnotesize
\begin{center}
%\begin{tabular}{|l|r|r|l|l|r|} \hline
\begin{tabular}{lrrlll} \hline
{Identifier} &
{$n$} &
 {$nnz(A)$} &
 {$normA$} &
 {$normb$} &
 {$cond2$} \\ 
% \multicolumn{6}{|c|}{\bf Well-conditioned problems} \\
\hline\Tstrut
HB/bcsstk27              &  1224 &   2.87$\times 10^4$ &   2.96$\times 10^7$ &   9.74$\times 10^5$   &     2.41$\times 10^4$ \\ 
Nasa/nasa2146            &  2146 &   3.72$\times 10^4$ &   2.79$\times 10^8$ &   9.05$\times 10^6$   &     1.72$\times 10^3$ \\ 
Cylshell/s1rmq4m1        &  5489 &   1.43$\times 10^5$ &   8.14$\times 10^6$ &   1.73$\times 10^5$   &     1.81$\times 10^6$ \\ 
MathWorks/Kuu            &  7102 &   1.74$\times 10^5$ &   4.73$\times 10^2$ &   5.01   &     1.58$\times 10^4$ \\ 
%UTEP/Dubcova1            & 16129 &   1.35$\times 10^5$ &   6.67$\times 10^{1}$ &   1.18   &     9.97$\times 10^2$ \\ 
Pothen/bodyy6            & 19366 &   7.71$\times 10^4$ &   1.09$\times 10^5$ &   9.81$\times 10^4$   &     9.91$\times 10^4$ \\ 
GHS$\_$psdef/wathen120   & 36441 &   3.01$\times 10^5$ &   1.52$\times 10^3$ &   2.66$\times 10^2$   &     9.58$\times 10^2$ \\ 
GHS$\_$psdef/jnlbrng1    & 40000 &   1.20$\times 10^5$ &   3.29$\times 10^{1}$ &   2.00$\times 10^{-1}$   &     1.83$\times 10^2$ \\ 
Williams/cant            & 62451 &   2.03$\times 10^6$ &   2.92$\times 10^5$ &   5.05$\times 10^3$   &     8.06$\times 10^3$ \\ 
UTEP/Dubcova2            & 65025 &   5.48$\times 10^5$ &   6.67$\times 10^{1}$ &   1.18   &     3.33 \\ 
Cunningham/qa8fm         & 66127 &   8.63$\times 10^5$ &   4.28$\times 10^{-3}$ &   9.51$\times 10^{-4}$   &     8.00 \\ 
Mulvey/finan512          & 74752 &   3.36$\times 10^5$ &   3.91$\times 10^2$ &   3.78$\times 10^{1}$   &     2.51$\times 10^{1}$ \\ 
GHS$\_$psdef/apache1     & 80800 &   3.11$\times 10^5$ &   8.10$\times 10^5$ &   6.76$\times 10^{-1}$   &     4.18$\times 10^2$ \\ 
Williams/consph          & 83334 &   3.05$\times 10^6$ &   6.61$\times 10^5$ &   7.20$\times 10^3$   &     1.25$\times 10^5$ \\ 
AMD/G2$\_$circuit        &150102 &   4.38$\times 10^5$ &   2.27$\times 10^4$ &   2.17$\times 10^4$   &     2.02$\times 10^4$ \\ 
%\multicolumn{6}{|c|}{\bf Ill-conditioned problems} \\
\hline\Tstrut
Boeing/msc01050          &  1050 &   1.51$\times 10^4$ &   2.58$\times 10^7$ &   1.90$\times 10^6$   &     4.58$\times 10^{15}$ \\ 
HB/bcsstk11              &  1473 &   1.79$\times 10^4$ &   1.21$\times 10^{10}$ &   7.05$\times 10^8$  &     2.21$\times 10^8$ \\ 
HB/bcsstk26              &  1922 &   1.61$\times 10^4$ &   1.68$\times 10^{11}$ &   8.99$\times 10^{10}$  &     1.66$\times 10^8$ \\ 
HB/bcsstk24              &  3562 &   8.17$\times 10^4$ &   5.28$\times 10^{14}$ &   4.21$\times 10^{13}$  &     1.95$\times 10^{11}$ \\ 
HB/bcsstk16              &  4884 &   1.48$\times 10^5$ &   4.12$\times 10^{10}$ &   9.22$\times 10^8$   &     4.94$\times 10^9$ \\ 
Cylshell/s2rmt3m1        &  5489 &   1.13$\times 10^5$ &   9.84$\times 10^5$ &   1.73$\times 10^4$ &     2.50$\times 10^8$ \\ 
Cylshell/s3rmt3m1        &  5489 &   1.13$\times 10^5$ &   1.01$\times 10^5$ &   1.73$\times 10^3$   &     2.48$\times 10^{10}$ \\ 
Boeing/bcsstk38          &  8032 &   1.82$\times 10^5$ &   4.50$\times 10^{11}$ &   4.04$\times 10^{11}$  &     5.52$\times 10^{16}$ \\ 
Boeing/msc10848          & 10848 &   6.20$\times 10^5$ &   4.58$\times 10^{13}$ &   6.19$\times 10^{11}$ &     9.97$\times 10^9$ \\ 
Oberwolfach/t2dah$\_$e   & 11445 &   9.38$\times 10^4$ &   2.20$\times 10^{-5}$ &   1.40$\times 10^{-5}$   &     7.23$\times 10^8$ \\ 
Boeing/ct20stif          & 52329 &   1.38$\times 10^6$ &   8.99$\times 10^{11}$ &   8.87$\times 10^{11}$   &     1.18$\times 10^{12}$ \\ 
DNVS/shipsec8            &114919 &   3.38$\times 10^6$ &   7.31$\times 10^{12}$ &   4.15$\times 10^{11}$  &     2.40$\times 10^{13}$ \\ 
Um/2cubes$\_$sphere      &101492 &   8.74$\times 10^5$ &   3.43$\times 10^{10}$ &   3.59$\times 10^{10}$  &     2.59$\times 10^8$ \\ 
GHS$\_$psdef/hood        &220542 &   5.49$\times 10^6$ &   2.23$\times 10^9$ &   1.51$\times 10^8$   &     5.35$\times 10^7$ \\ 
Um/offshore              &259789 &   2.25$\times 10^6$ &   1.44$\times 10^{15}$ &   1.16$\times 10^{15}$   &     4.26$\times 10^9$ \\ 
\hline
\end{tabular}
\end{center}
}
\end{table}

We use Algorithm~\ref{alg:gmres-ir-ic} to explore whether
we can recover (close to) double precision accuracy using
preconditioners computed in fp16 arithmetic, although we are 
aware that in practice much less accuracy in the computed solution may be sufficient
(indeed, in many practical situations, inaccuracies in the supplied
data may mean low precision accuracy in the solution is all that can be justified). We thus use
two precisions: $u_{\ell}=u_{16}$  for the factorization and $u_r = u_g = u_p = u_{64}$.
Algorithm~\ref{alg:gmres-ir-ic} terminates when the
normwise backward error for the computed solution satisfies
\begin{equation}
\label{eq:resid}
res = \frac{\| b - A  x\|_\infty}{\| A\|_\infty \|  x\|_\infty + \| b\|_\infty} \le \delta.
 \end{equation}
In our experiments, we set $\delta = 10^3 \times u_{64}$.
The implementations of CG and GMRES  used are {\tt MI21} and {\tt MI24}, respectively,
from the HSL software library \cite{hsl:2023}.
Except for the results in Table~\ref{T:changingtolerance}, 
the CG and the GMRES convergence tolerance is $\delta_{krylov} = u_{64}^{1/4}$
and the limit on the number of iterations for each application of CG and GMRES is 1000.

The numerical experiments are performed on a Windows 11-Pro-based machine with 
an Intel(R) Core(TM) i5-10505 CPU processor (3.20GHz).
Our results are for the level-based incomplete Cholesky factorization $IC(\ell)$ for a range of values of $\ell \ge 0$.
The sparsity pattern of $L$ is computed using the approach of Hysom and Pothen~\cite{hypo:98}.
This computes the pattern of each row of $L$ independently. % and thus can be done in parallel.
Our software is written in Fortran and compiled using the NAG compiler (Version 7.1, Build 7118).
This is currently the only Fortran compiler that supports the use of fp16.
The NAG documentation states that their half precision implementation conforms to the IEEE standard. In addition,
using the {\tt -round$\_$hreal} option, all half-precision operations are 
rounded to half precision, both at compile time and runtime.
{Section 3 of \cite{hipr:2019} presents an insightful discussion on rounding every operation or rounding kernels (see also \cite{hima:2022})}.
Because of all the conversions needed, half precision is slower using the NAG compiler
than single precision and so timings are not useful.
%Note that, as far as we are aware, current C and C++ compilers also only simulate half precision arithmetic.

We refer to the $IC(\ell)$ factorizations computed using half and double precision arithmetic as
fp16-$IC(\ell)$ and fp64-$IC(\ell)$, respectively.
The key difference between the fp16 and fp64 versions of our $IC(\ell)$ software is
that for the former, during the incomplete factorization,
we incorporate the safe checks for the scaling and update operations (as discussed
in Section~\ref{sec:breakdowns}); for the fp64 version,
only tests for B1 breakdowns are performed (B2 and B3 breakdowns
are not encountered in our double precision experiments). In addition, the fp16 version allows the preconditioner to be applied
in either half or double precision arithmetic; the former is for IC-LU-IR and the latter for IC-Krylov-IR.
In the IC-Krylov-IR case,  {preconditioning
is performed in double precision. 
For $L$ computed using fp16, there are
two straightforward ways of handling solving systems with $L$ and $L^T$ in double precision.}
The first makes an explicit copy of $L$ by casting the data into double precision
but this negates the important benefit that half precision offers of reducing memory requirements.
The alternative casts the entries on the fly. This
is straightforward to incorporate into a serial triangular solve routine, and 
only requires a temporary double precision array of length $n$.

In the results tables, $nnz(L)$ is the number of entries in the incomplete factor $L$;
$iouter$ denotes the number of iterative refinement steps (that is, the number of times the
loop starting at Line 3 in Algorithm~\ref{alg:gmres-ir-ic} is executed)
and $totits$ is the total number of CG (or GMRES) iterations performed;
$resint$ is  (\ref{eq:resid}) with $x =  L^{-T}L^{-1}b$
(because $L$ depends on the precision used
to compute it, $resint$ is precision dependent),
and $resfinal$ is (\ref{eq:resid}) for the final computed solution;
$nmod$ and $nofl$ are the numbers of times problems B1 and B3 occur during the incomplete factorization, with
the latter for fp16 only. 
Increasing the global shift and restarting the factorization after the detection of problem B1
avoided problem B2 in all our tests.
Problem B3 can occur in column $k$ if after the shift, 
the diagonal entry $l_{kk}$ is close to the shift $\alpha$ 
and $|l_{ik}|/x_{max} \ge \alpha $ for some $i > k$.
In all the experiments on well-conditioned problems, we found $nofl = 0$ and so this statistic is
not included in the corresponding tables of results.
As expected, $nmod >0$ can occur for fp16 and fp64, and for both well-conditioned
and ill-conditioned examples.

\subsection{Results for IC-LU-IR}

IC-LU-IR is attractive because if fp16 arithmetic is used then the application of the preconditioner is performed in half precision arithmetic.
Table~\ref{T:LU-IRw} reports results for IC-LU-IR with the $IC(3)$
preconditioner.
The iteration count is limited to 1000. When using fp16 and fp64 arithmetic, we were unable to 
to achieve the requested accuracy
within this limit for some problems. Additionally, for problem  Williams/consph, the refinement procedure 
diverges and the process is stopped when 
the norm of the residual approaches $x_{max}$ in double precision.
For the ill-conditioned problems, results are given only 
the ones that were successfully solved;
for the other test examples, we failed to achieve convergence.

\begin{table}[htbp]
\caption{IC-LU-IR results using the $IC(3)$ preconditioner.
$resint$ and $resfinal$ are the initial and final scaled residuals;
$nnz(L)$ is the number of entries in the $IC(3)$ factor; $iters$ is the number of refinement steps; and
$nmod$ denotes the number of times problem B1 occurs during the factorization
(for fp64-$IC(3)$ it is equal to 0 for all our test cases and is omitted).
$> 1000$ indicates the requested accuracy was not obtained within the  iteration limit.
$\dag$ indicates the refinement procedure breaks down.
}
\label{T:LU-IRw}\vspace{3mm}
{
\footnotesize
\begin{center}
%\begin{tabular}{|l|r|r|r|r|r|} \hline
\begin{tabular}{lrlrrr} \hline
\multicolumn{6}{c}{Preconditioner fp16-$IC(3)$} \\
{Identifier} &
{$resinit$} &
{$resfinal$} &
{$nnz(L)$} &
{$iters$} &
{$nmod$} \\
\hline\Tstrut
HB/bcsstk27              &     8.27$\times 10^{-5}$  &      6.11$\times 10^{-14}$  &      4.88$\times 10^4$ &    13      &     0  \\ 
Nasa/nasa2146            &     9.07$\times 10^{-5}$  &      1.89$\times 10^{-13}$  &      7.89$\times 10^4$ &    15      &     0  \\ 
Cylshell/s1rmq4m1        &     5.81$\times 10^{-5}$  &      6.77$\times 10^{-9}$  &     3.15$\times 10^5$ &     $>1000$      &     0  \\ 
MathWorks/Kuu            &     1.74$\times 10^{-4}$  &      2.17$\times 10^{-13}$  &      7.65$\times 10^5$ &   194      &     0  \\ 
Pothen/bodyy6            &     7.39$\times 10^{-3}$  &      2.22$\times 10^{-13}$  &      1.76$\times 10^5$ &   822      &     2  \\ 
GHS$\_$psdef/wathen120   &     4.34$\times 10^{-4}$  &      1.72$\times 10^{-14}$  &      8.30$\times 10^5$ &     5      &     0  \\ 
GHS$\_$psdef/jnlbrng1    &     1.41$\times 10^{-3}$  &      6.48$\times 10^{-14}$  &      2.77$\times 10^5$ &    16      &     0  \\ 
Williams/cant            &     3.26$\times 10^{-4}$  &      4.61$\times 10^{-8}$  &      9.95$\times 10^6$ &     $>1000$      &     0  \\ 
UTEP/Dubcova2            &     1.54$\times 10^{-3}$  &      2.19$\times 10^{-13}$  &      6.22$\times 10^6$ &   736      &     0  \\ 
Cunningham/qa8fm         &     3.54$\times 10^{-4}$  &      3.42$\times 10^{-15}$  &      5.14$\times 10^6$ &     5      &     0  \\ 
Mulvey/finan512          &     3.50$\times 10^{-4}$  &      2.59$\times 10^{-15}$  &      4.08$\times 10^6$ &     5      &     0  \\ 
GHS$\_$psdef/apache1     &     8.52$\times 10^{-5}$  &      1.05$\times 10^{-7}$  &      1.54$\times 10^6$ &   $>1000$     &     0  \\ 
Williams/consph          &     1.05$\times 10^{-4}$  &          $~~~~~\dag$  &      2.02$\times 10^7$ &     $\dag~~$   & 0    \\ 
AMD/G2$\_$circuit        &     8.16$\times 10^{-4}$  &      7.71$\times 10^{-8}$  &      1.04$\times 10^6$ &     $>1000$     &     0  \\ 
\hline\Tstrut
Oberwolfach/t2dah$\_$e      &     6.81$\times 10^{-4}$  &      1.42$\times 10^{-14}$  &      3.29$\times 10^5$ &     5  &     0     \\ 
Um/2cubes$\_$sphere         &     1.02$\times 10^{-3}$  &      9.10$\times 10^{-16}$  &      8.70$\times 10^6$ &     5  &     0     \\ 
\hline\hline
\multicolumn{6}{c}{Preconditioner fp64-$IC(3)$} \\
{Identifier} &
{$resinit$} &
{$resfinal$} &
{$nnz(L)$} &
{$iters$} &
{} \\
\cline{1-5}\Tstrut
HB/bcsstk27              &     3.93$\times 10^{-5}$  &      2.06$\times 10^{-14}$  &      4.88$\times 10^4$ &     7      &  \\ 
Nasa/nasa2146            &     4.32$\times 10^{-5}$  &      6.29$\times 10^{-14}$  &      7.89$\times 10^4$ &    15      &       \\ 
Cylshell/s1rmq4m1        &     3.55$\times 10^{-5}$  &      4.11$\times 10^{-9}$  &     3.15$\times 10^5$ &     $>1000$      &      \\ 
MathWorks/Kuu            &     2.50$\times 10^{-4}$  &      1.95$\times 10^{-13}$  &      7.64$\times 10^5$ &   110      &       \\ 
Pothen/bodyy6            &     9.31$\times 10^{-3}$  &      2.20$\times 10^{-13}$  &      1.76$\times 10^5$ &   $>1000$  &       \\ 
GHS$\_$psdef/wathen120   &     1.77$\times 10^{-4}$  &      5.62$\times 10^{-15}$  &      8.30$\times 10^5$ &     5      &       \\ 
GHS$\_$psdef/jnlbrng1    &     9.40$\times 10^{-3}$  &      1.14$\times 10^{-13}$  &      2.77$\times 10^5$ &    14      &      \\ 
Williams/cant            &     3.19$\times 10^{-4}$  &      2.55$\times 10^{-8}$  &      9.95$\times 10^6$ &     $>1000$      &      \\ 
UTEP/Dubcova2            &     1.65$\times 10^{-3}$  &      2.17$\times 10^{-13}$  &      6.22$\times 10^6$ &   709      &      \\ 
Cunningham/qa8fm         &     3.29$\times 10^{-4}$  &      1.49$\times 10^{-15}$  &      5.14$\times 10^6$ &     4      &      \\ 
Mulvey/finan512          &     3.72$\times 10^{-4}$  &      9.22$\times 10^{-15}$  &      4.08$\times 10^6$ &     4      &       \\ 
GHS$\_$psdef/apache1     &     8.28$\times 10^{-5}$  &      2.83$\times 10^{-8}$  &      1.54$\times 10^6$ &   $>1000$     &       \\ 
Williams/consph          &     3.38$\times 10^{-4}$  &          $~~~~~\dag$  &      2.02$\times 10^7$ &     $\dag~~$   &    \\ 
AMD/G2$\_$circuit        &     8.42$\times 10^{-5}$  &      2.64$\times 10^{-8}$  &      1.04$\times 10^6$ &     $>1000$     &       \\ 
\cline{1-5}\Tstrut
HB/bcsstk24                 &     4.01$\times 10^{-7}$  &      2.21$\times 10^{-13}$  &      2.77$\times 10^5$ &   801    &        \\ 
HB/bcsstk16                 &     7.01$\times 10^{-4}$  &      1.81$\times 10^{-9} $  &      4.89$\times 10^6$ &  $>1000$ &         \\ 
Oberwolfach/t2dah$\_$e      &     5.49$\times 10^{-6}$  &      5.95$\times 10^{-14}$  &      3.29$\times 10^5$ &     4  &         \\ 
Um/2cubes$\_$sphere         &     2.25$\times 10^{-5}$  &      1.10$\times 10^{-13}$  &      8.70$\times 10^6$ &     3  &         \\ 
\cline{1-5}
\end{tabular}
\end{center}
}
\end{table}

\subsection{Dependence of the iteration counts on the CG tolerance}
The results in Table~\ref{T:changingtolerance} illustrate the dependence of the number of iterative
refinement steps and the total iteration count on the
convergence tolerance $\delta_{krylov}$ used by CG within IC-CG-IR. 
The preconditioner $IC(3)$ computed using fp16 arithmetic. 
For three examples, we test
$\delta_{krylov}$ ranging from $10^{-10}$ to $10^{-1}$. The first problem is well-conditioned
and the other two are ill-conditioned. As expected, the number of outer iterations increases slowly
with $\delta_{krylov}$, but the precise choice of $\delta_{krylov}$ is not critical. 
The results  confirm the choice of $u_{64}^{1/4}$,
which is used for all remaining experiments.
Similar results are obtained for IC-GMRES-IR, emphasising that an advantage of the IR 
approach is that it can  significantly reduce the maximum number of Krylov vectors used.
\begin{table}[htbp]
\caption{The effects of changing the CG convergence tolerance $\delta_{krylov}$  used in IC-CG-IR.
The preconditioner is fp16-$IC(3)$.
$iouter$ and $totits$ denote the number of outer iterations and the total number of CG iterations, respectively.
}
\label{T:changingtolerance}\vspace{3mm}
{
\footnotesize
\begin{center}
\begin{tabular}{|l|l|l|l|l|l|l|l|l|l|l|} \hline
\hline\Tstrut
   &    \multicolumn{10}{|c|} {\tt UTEP/Dubcova2} \\
$\delta_{krylov}$  &  $10^{-10}$& $10^{-9}$ & $10^{-8}$ & $10^{-7}$ & $10^{-6}$ & $10^{-5}$ & $10^{-4}$ & $10^{-3}$ & $10^{-2}$ & $10^{-1}$   \\
$iouter$           & 2 & 2        & 2 & 2 &  2 &   3 &  3 &  4 &  6 &  8  \\ 
$totits$           & 79 & 73 & 64 & 58 & 49 &  68 & 55 & 57 & 58 &  65  \\ 
\hline
\hline\Tstrut
   &    \multicolumn{10}{|c|} {\tt HB/bcsstk26} \\
$\delta_{krylov}$ &  $10^{-10}$& $10^{-9}$ & $10^{-8}$ & $10^{-7}$ & $10^{-6}$ & $10^{-5}$ & $10^{-4}$ & $10^{-3}$ & $10^{-2}$ & $10^{-1}$   \\
$iouter$          &  2 & 2 &    2 & 2    & 2 &    3 &    3 &   4 & 6 &   9 \\ 
$totits$          &  134 & 121 & 107 & 93 & 78 & 103 & 81 & 84 & 97 & 93  \\ 
\hline
\hline\Tstrut
   &    \multicolumn{10}{|c|} {\tt Cylshell/s2rmt3m1} \\
$\delta_{krylov}$  &  $10^{-10}$& $10^{-9}$ & $10^{-8}$ & $10^{-7}$ & $10^{-6}$ & $10^{-5}$ & $10^{-4}$ & $10^{-3}$ & $10^{-2}$ & $10^{-1}$   \\
$iouter$           &  2 &       2 &    2 &     2 &    2 & 3 &   3 & 4 & 5 & 8  \\ 
$totits$           &  132 & 130 & 125 & 120 & 116 & 94 &   83 & 117 &  85 & 150 \\ 
\hline
\hline\Tstrut
   &    \multicolumn{10}{|c|} {\tt GHS$\_$psdef/wathen120} \\
$\delta_{krylov}$  &  $10^{-10}$& $10^{-9}$ & $10^{-8}$ & $10^{-7}$ & $10^{-6}$ & $10^{-5}$ & $10^{-4}$ & $10^{-3}$ & $10^{-2}$ & $10^{-1}$   \\
$iouter$           &  2 &   2 &  2 &     2 &  2 & 2 &     3 &    3 &    6 &    6 \\
$totits$           &  9 &   8 &  7 &     6 &  6 &  5 &    6 &    5 &    6 &    6 \\
\hline
\hline\Tstrut
   &    \multicolumn{10}{|c|} {\tt GHS$\_$psdef/jnlbrng1} \\
$\delta_{krylov}$  &  $10^{-10}$& $10^{-9}$ & $10^{-8}$ & $10^{-7}$ & $10^{-6}$ & $10^{-5}$ & $10^{-4}$ & $10^{-3}$ & $10^{-2}$ & $10^{-1}$   \\
$iouter$           &  2  &  2 &  2 &     2 &  2 &  2 &    3 &     4 &   5 &    7 \\
$totits$           &  20 &  18 & 16 &   14 &  12 & 15 &   12 &   12 &   11 &  12 \\
\hline
\hline\Tstrut
   &    \multicolumn{10}{|c|} {\tt Oberwolfach/t2dah$\_$e} \\
$\delta_{krylov}$  &  $10^{-10}$& $10^{-9}$ & $10^{-8}$ & $10^{-7}$ & $10^{-6}$ & $10^{-5}$ & $10^{-4}$ & $10^{-3}$ & $10^{-2}$ & $10^{-1}$   \\
$iouter$           &   2 &   2 &  2 &  2 &  2 &   3 &   3 &       3 &    5 &   5 \\
$totits$           &   9 &   8 &  7 &  6 &  5 &   7 &   6 &       5 &    5 &   5 \\
\hline
\end{tabular}
\end{center}

}
\end{table}

\subsection{Results for $IC(0)$}
\begin{table}[htbp]
\caption{Results for IC-CG-IR using an $IC(0)$ preconditioner: well-conditioned problems.
$resint$ and $resfinal$ are the initial and final scaled residuals.
$nnz(L)$ is the number of entries in the $IC(0)$ factor. $iouter$ and 
$totits$ denote the number of outer iterations and the total number of CG iterations, respectively
$>1000$ indicates CG tolerance not reached on outer iteration $iouter$.
$nmod$ denotes the number of times problem B1 occurs during the factorization.
A count in bold indicates the fp16 result is within 10 per cent of (or is better than) the corresponding fp64 result.
}

\label{T:IC(0)+CGw}\vspace{3mm}
{
\footnotesize
\begin{center}
%\begin{tabular}{|l|r|r|r|r|r|r|r|} \hline
\begin{tabular}{lrrrrrr} \hline
\multicolumn{7}{c}{Preconditioner fp16-$IC(0)$} \\
{Identifier} &
{$resinit$} &
{$resfinal$} &
{$nnz(L)$} &
 {$iouter$} &
 {$totits$} &
 {$nmod$} \\
\hline\Tstrut
HB/bcsstk27              &     4.19$\times 10^{-4}$  &      2.02$\times 10^{-14}$  &      2.87$\times 10^4$ &     3  &    {\bf 35}    &     0  \\ 
Nasa/nasa2146            &     3.26$\times 10^{-4}$  &      3.82$\times 10^{-15}$  &      3.72$\times 10^4$ &     3  &    {\bf 24}    &     0  \\ 
Cylshell/s1rmq4m1        &     1.49$\times 10^{-4}$  &      2.16$\times 10^{-14}$  &      1.15$\times 10^5$ &     3  &   210    &     0  \\ 
MathWorks/Kuu            &     3.27$\times 10^{-3}$  &      1.38$\times 10^{-14}$  &      1.43$\times 10^5$ &     3  &   274    &     4  \\ 
Pothen/bodyy6            &     1.71$\times 10^{-4}$  &      1.58$\times 10^{-16}$  &      7.03$\times 10^4$ &     4  &   178    &     2  \\ 
GHS$\_$psdef/wathen120   &     1.09$\times 10^{-2}$  &      1.09$\times 10^{-13}$  &      3.01$\times 10^5$ &     3  &    {\bf 17}    &     0  \\ 
GHS$\_$psdef/jnlbrng1    &     9.98$\times 10^{-3}$  &      1.18$\times 10^{-14}$  &      1.20$\times 10^5$ &     3  &    {\bf 40}    &     0  \\ 
Williams/cant            &     4.62$\times 10^{-3}$  &      2.81$\times 10^{-8}$   &      1.46$\times 10^6$ &     2  &   $>1000$    &     9  \\ 
UTEP/Dubcova2            &     5.55$\times 10^{-3}$  &      1.86$\times 10^{-13}$  &      4.19$\times 10^5$ &     3  &   {\bf 225}    &     0  \\ 
Cunningham/qa8fm         &     3.59$\times 10^{-3}$  &      3.52$\times 10^{-16}$  &      8.63$\times 10^5$ &     4  &    {\bf 14}    &     0  \\ 
Mulvey/finan512          &     2.88$\times 10^{-3}$  &      1.25$\times 10^{-14}$  &      3.36$\times 10^5$ &     3  &    {\bf 14}    &     0  \\ 
GHS$\_$psdef/apache1     &     3.61$\times 10^{-4}$  &      2.54$\times 10^{-14}$  &      3.11$\times 10^5$ &     2  &   { 274}    &     0  \\ 
Williams/consph          &     8.59$\times 10^{-5}$  &      1.62$\times 10^{-14}$  &      3.05$\times 10^6$ &     3  &   {\bf 435}    &     7  \\ 
AMD/G2$\_$circuit        &     7.75$\times 10^{-4}$  &      9.01$\times 10^{-16}$  &      4.38$\times 10^5$ &     4  &   {\bf 842}     &     0 \\ 
\hline\hline
\multicolumn{7}{c}{Preconditioner fp64-$IC(0)$} \\
{Identifier} &
{$resinit$} &
{$resfinal$} &
{$nnz(L)$} &
 {$iouter$} &
 {$totits$} &
 {$nmod$}  \\
\hline\Tstrut
HB/bcsstk27              &     4.31$\times 10^{-4}$  &      2.42$\times 10^{-14}$  &      2.87$\times 10^4$ &     3  &    35    &     0       \\ 
Nasa/nasa2146            &     3.06$\times 10^{-4}$  &      4.74$\times 10^{-15}$  &      3.72$\times 10^4$ &     3  &    24    &     0       \\ 
Cylshell/s1rmq4m1        &     1.61$\times 10^{-4}$  &      1.57$\times 10^{-14}$  &      1.43$\times 10^5$ &     3  &   176    &     0       \\ 
MathWorks/Kuu            &     9.47$\times 10^{-4}$  &      2.47$\times 10^{-14}$  &      1.74$\times 10^5$ &     3  &   130    &     0       \\ 
Pothen/bodyy6            &     8.68$\times 10^{-5}$  &      1.40$\times 10^{-16}$  &      7.71$\times 10^4$ &     4  &   138    &     0       \\ 
GHS$\_$psdef/wathen120   &     1.09$\times 10^{-2}$  &      1.10$\times 10^{-13}$  &      3.01$\times 10^5$ &     3  &    17    &     0     \\ 
GHS$\_$psdef/jnlbrng1    &     9.96$\times 10^{-3}$  &      1.05$\times 10^{-14}$  &      1.20$\times 10^5$ &     3  &    40    &     0     \\ 
Williams/cant            &     2.63$\times 10^{-3}$  &      4.82$\times 10^{-8}$   &      2.03$\times 10^6$ &     2  &  $>1000$    &     8       \\ 
UTEP/Dubcova2            &     5.36$\times 10^{-3}$  &      1.40$\times 10^{-13}$  &      5.48$\times 10^5$ &     3  &   227    &     0       \\ 
Cunningham/qa8fm         &     3.57$\times 10^{-3}$  &      1.23$\times 10^{-16}$  &      8.63$\times 10^5$ &     4  &    14    &     0       \\ 
Mulvey/finan512          &     2.92$\times 10^{-3}$  &      4.18$\times 10^{-14}$  &      3.36$\times 10^5$ &     3  &    14    &     0      \\ 
GHS$\_$psdef/apache1     &     3.63$\times 10^{-4}$  &      1.81$\times 10^{-14}$  &      3.11$\times 10^5$ &     2  &   244    &     0    \\ 
Williams/consph          &     8.26$\times 10^{-5}$  &      2.46$\times 10^{-14}$  &      3.05$\times 10^6$ &     3  &   440    &     7       \\ 
AMD/G2$\_$circuit        &     5.44$\times 10^{-4}$  &      6.55$\times 10^{-16}$  &      4.38$\times 10^5$ &     4  &   855    &     0   \\ 
\hline
\end{tabular}
\end{center}
}
\end{table}

\begin{table}[htbp]
\caption{Results for IC-CG-IR and IC-GMRES-IR using an $IC(0)$ preconditioner: ill-conditioned problems.
$resint$ is the initial scaled residual;
$resfinal$ is the final IC-CG-IR scaled residual.
$nnz(L)$ is the number of entries in the $IC(0)$ factor. $iouter$ and 
$totits$ denote the number of outer iterations and the total number of CG iterations with the GMRES
statistics in parentheses.
$>1000$ indicates CG (or GMRES) tolerance not reached on outer iteration $iouter$.
$nmod$ denotes the number of times problem B1 occurs during the factorization.
A count in bold indicates the fp16 result is within 10 per cent of (or is better than) the corresponding fp64 result.
$^\ast$ denotes early termination of CG.}
\label{T:IC(0)+CGi}\vspace{3mm}
{
\footnotesize
\begin{center}
\begin{tabular}{lrlrrlrlr} \hline
\multicolumn{9}{c}{Preconditioner fp16-$IC(0)$} \\
\multicolumn{1}{c}{Identifier} &
\multicolumn{1}{c}{$resinit$} &
\multicolumn{1}{c}{$resfinal$} &
\multicolumn{1}{c}{$nnz(L)$} &
\multicolumn{2}{c}{$iouter$} &
\multicolumn{2}{c}{$totits$} &
\multicolumn{1}{c}{$nmod$} \\
\hline\Tstrut
Boeing/msc01050          &     1.45$\times 10^{-5}$  &      1.20$\times 10^{-13}$  &      9.49$\times 10^3$ &     3  & (3)&  572 & ({\bf 333})    &     7  \\ 
HB/bcsstk11              &     5.37$\times 10^{-4}$  &      1.66$\times 10^{-13}$  &      1.52$\times 10^4$ &     3  & (3)&   914 &(644)    &     5  \\ 
HB/bcsstk26              &     6.70$\times 10^{-4}$  &      1.18$\times 10^{-16}$  &      1.33$\times 10^4$ &     4  & (4)&  {\bf 488} &({\bf 422})   &     2  \\ 
HB/bcsstk24              &     4.50$\times 10^{-5}$  &      1.41$\times 10^{-12}$  &      7.97$\times 10^4$ &     3  & (3)&  $> 1000$ & ({\bf 834})  &     3  \\ 
HB/bcsstk16              &     3.36$\times 10^{-3}$  &      3.50$\times 10^{-14}$  &      1.27$\times 10^5$ &     3  & (3)&  88 & (80)  &     4  \\ 
Cylshell/s2rmt3m1        &     1.41$\times 10^{-4}$  &      9.88$\times 10^{-15}$  &      1.02$\times 10^5$ &     3  & (3)&   { 347} &({\bf 473})    &     0  \\ 
Cylshell/s3rmt3m1        &     1.13$\times 10^{-5}$  &      1.64$\times 10^{-14}$  &      1.02$\times 10^5$ &     3  & (3)&  $> 1000$ & ($> 1000$)  &     3  \\ 
Boeing/bcsstk38          &     3.27$\times 10^{-2}$  &      1.58$\times 10^{-10}$  &      1.63$\times 10^5$ &     3& (3)  &  $> 1000$ &($> 1000$)    &     8  \\ 
Boeing/msc10848          &     5.46$\times 10^{-6}$  &      1.96$\times 10^{-14}$  &      6.18$\times 10^5$ &     3  & (3)&  $> 1000$ &({\bf 622})   &     2  \\ 
Oberwolfach/t2dah$\_$e   &     4.87$\times 10^{-2}$  &      5.61$\times 10^{-9}$$^\ast$  &      9.38$\times 10^4$ &     1  & (4)&    {\bf 15} & ({\bf 34})  &     0  \\ 
Boeing/ct20stif          &     3.82$\times 10^{-3}$  &      1.83$\times 10^{-9}$   &      1.30$\times 10^6$ &     3  & (3)&  $> 1000$ &($> 1000$)   &     7  \\ 
DNVS/shipsec8            &     2.31$\times 10^{-3}$  &      2.40$\times 10^{-9}$   &      1.53$\times 10^6$ &     3  & (3)&  $> 1000$ & ($> 1000$)  &     8  \\ 
Um/2cubes$\_$sphere      &     1.88$\times 10^{-2}$  &      1.77$\times 10^{-14}$  &      8.74$\times 10^5$ &     3  & (3)&    {\bf 13} &({\bf 11})  &     0  \\ 
GHS$\_$psdef/hood        &     1.93$\times 10^{-3}$  &      5.01$\times 10^{-17}$  &      5.06$\times 10^6$ &     4  & (3)&  { 620} &({\bf 346})   &     2  \\ 
Um/offshore              &     1.49$\times 10^{-2}$  &      2.88$\times 10^{-16}$  &      2.25$\times 10^6$ &     4  & (3)&   {\bf 602} &({\bf 101})   &     0  \\ 
\hline
\multicolumn{9}{c}{Preconditioner fp64-$IC(0)$} \\
\multicolumn{1}{c}{Identifier} &
\multicolumn{1}{c}{$resinit$} &
\multicolumn{1}{c}{$resfinal$} &
\multicolumn{1}{c}{$nnz(L)$} &
\multicolumn{2}{c}{$iouter$} &
\multicolumn{2}{c}{$totits$} &
\multicolumn{1}{c}{$nmod$} \\
\hline\Tstrut
Boeing/msc01050          &     6.03$\times 10^{-3}$  &      3.11$\times 10^{-14}$  &      1.51$\times 10^4$ &     3  & (3)&   555  & (485) &     8       \\ 
HB/bcsstk11              &     7.49$\times 10^{-4}$  &      7.44$\times 10^{-14}$  &      1.79$\times 10^4$ &     3  & (3)&  902  & (475)  &     4       \\ 
HB/bcsstk26              &     1.59$\times 10^{-3}$  &      2.37$\times 10^{-16}$  &      1.61$\times 10^4$ &     4  & (4)&   592  &  (504) &     4       \\ 
HB/bcsstk24              &     4.59$\times 10^{-5}$  &      1.11$\times 10^{-13}$  &      8.17$\times 10^4$ &     3  & (3)&  $> 1000$  & (812)  &     3       \\ 
HB/bcsstk16              &     2.76$\times 10^{-3}$  &      2.66$\times 10^{-14}$  &      1.48$\times 10^5$ &     3  & (3)&    68  & (66)  &     0       \\ 
Cylshell/s2rmt3m1        &     1.45$\times 10^{-4}$  &      1.10$\times 10^{-14}$  &      1.13$\times 10^5$ &     3  & (3)&   314  & (455)  &     0       \\ 
Cylshell/s3rmt3m1        &     1.14$\times 10^{-5}$  &      7.88$\times 10^{-15}$  &      1.13$\times 10^5$ &     3  & (3)&   584  & (719) &     0       \\ 
Boeing/bcsstk38          &     1.21$\times 10^{-2}$  &      1.69$\times 10^{-13}$  &      1.82$\times 10^5$ &     3  & (3)&  $> 1000$  & ($> 1000$)  &     7       \\ 
Boeing/msc10848          &     1.07$\times 10^{-5}$  &      1.38$\times 10^{-15}$  &      6.20$\times 10^5$ &     3  & (3)&  $> 1000$  & (760)  &     3       \\ 
Oberwolfach/t2dah$\_$e   &     4.86$\times 10^{-2}$  &      1.52$\times 10^{-16}$$^\ast$ &  9.38$\times 10^4$ &     1  & (4)&    15 & (34)   &     0    \\ 
Boeing/ct20stif          &     3.65$\times 10^{-3}$  &      3.20$\times 10^{-12}$  &      1.38$\times 10^6$ &     3  & (3)&  $> 1000$ & ($> 1000$)   &     7       \\ 
DNVS/shipsec8            &     2.75$\times 10^{-4}$  &      6.00$\times 10^{-13}$  &      3.38$\times 10^6$ &     3  & (3)&  $> 1000$ & ($> 1000$)   &     4       \\ 
Um/2cubes$\_$sphere      &     1.89$\times 10^{-2}$  &      2.55$\times 10^{-14}$  &      8.74$\times 10^5$ &     3  & (3)&    13  & (11)  &     0    \\ 
GHS$\_$psdef/hood        &     9.82$\times 10^{-4}$  &      1.73$\times 10^{-13}$  &      5.49$\times 10^6$ &     3  & (3)&   546  & (377)  &     2    \\ 
Um/offshore              &     1.51$\times 10^{-2}$  &      1.92$\times 10^{-13}$  &      2.25$\times 10^6$ &     4  & (3)&   590  & (100)  &     0      \\ 
\hline
\end{tabular}
\end{center}
}
\end{table}
As already remarked, $IC(0)$ is a very simple preconditioner but one that is frequently reported
on in publications. Results for IC-CG-IR using an $IC(0)$ preconditioner
computed in half and double precision arithmetics are given in Tables~\ref{T:IC(0)+CGw} and \ref{T:IC(0)+CGi}
for well-conditioned and ill-conditioned problems, respectively. If the total iteration count ($totits$) for fp16
is within 10 per cent of the count for fp64 (or is less than the fp64 count) then it is highlighted in bold. {Note that the number of
entries in $L$ is equal to the number of entries 
in the matrix that is being
factorized, that is, $nnz(L)= nnz(\widehat A^{(l)})$ for fp16
and $nnz(L) = nnz(A)$ for fp64. The difference between them is the number of entries that underflow and are dropped
when the scaled matrix is squeezed into half precision. For the well-conditioned problems, 
the only problem for which $\widehat A^{(l)}$ is significantly
sparser than $A$ is Williams/cant but for some of the ill-conditioned
problems (including Boeing/msc01050 and DNVS/shipsec8), $nnz(L)$ is much smaller for fp16 than for
fp64.}

From the tables we see that, for well-conditioned problems, using fp16 arithmetic to compute the $IC(0)$ factorization
is often as good as using fp64 arithmetic. For problem Williams/cant,  on the second 
outer iteration the CG method
fails to converge within 1000 iterations for both the fp16 and the fp64) preconditioners.
$nofl$ is omitted from Tables~\ref{T:IC(0)+CGw} and \ref{T:IC(0)+CGi} because it was 0 for all our test examples.
However, for many problems (particularly the ill-conditioned ones), $nmod > 0$
for both half precision and double precision and this can lead to a poor quality 
preconditioner, indicated by high iteration counts, with the limit of 1000
iterations being exceeded on the third outer iteration
for a number of test examples (such as Boeing/bcsstk38 and Boeing/msc10848). 
Although  the requested accuracy is not achieved, there are still
significant reductions in the initial residual so the preconditioners may be acceptable
if less accuracy is required. Nevertheless, the fp16 performance is often
competitive with that of fp64.
For problem Oberwolfach/t2dah$\_$e, the CG  algorithm 
terminates before the requested accuracy has been achieved; this is because the curvature encountered 
within the CG algorithm is found to be too small, triggering an error return.

We have also tested IC-GMRES-IR with the $IC(0)$ preconditioners.
For the well-conditioned problems, the iteration counts using GMRES are similar to those for CG.
For the ill-conditioned examples, the IC-GMRES-IR counts are given in parentheses in 
the $iouter$ and $totits$ columns of Table~\ref{T:IC(0)+CGi}.
Our findings are broadly consistent  with those reported in \cite{hipr:2021}. Although
IC-GMRES-IR  can require fewer iterations than IC-CG-IR,
it is important to remember that a GMRES iteration is more expensive than a CG iteration
so simply comparing  counts may be misleading.

In Table~\ref{T:IC(0)+GMRESi_pure}, we give results for $itmax = 1$, that is,
preconditioned GMRES is not applied within a refinement loop (we denote this by IC-GMRES).  
The initial residual, the number of entries in the factor, and the number of 
occurrences of problem B1  are as in Table~\ref{T:IC(0)+CGi} and so are not included.
We also report  the maximum number $maxbasis$ of
GMRES iterations performed on an outer iteration of IC-GMRES-IR
using the same fp16-$IC(0)$ preconditioner.
We see that the total iteration count for IC-GMRES-IR
can exceed the count for IC-GMRES but, for many of the problems
(including Boeing/msc01050, HB/bcsstk26 and Oberwolfach/t2dah$\_$e),
the maximum size of the constructed Krylov basis is 
smaller for IC-GMRES-IR than for IC-GMRES. Another option would be to use 
restarted GMRES, which would involve selecting an appropriate restart parameter.
We do not consider this here.

\begin{table}[htbp]
\caption{Results for IC-GMRES (without outer steps) using the fp16-$IC(0)$ preconditioner.
$resfinal$ is the final  scaled residual.
$totits$ denotes the number of GMRES iterations, with the counts for IC-GMRES-IR in parentheses.
$>2000$ indicates convergence is not achieved within $2000$ iterations.
$maxbasis$ is the maximum number of
GMRES iterations performed on an outer iteration of IC-GMRES-IR.
}
\label{T:IC(0)+GMRESi_pure}\vspace{3mm}
{
\footnotesize
\begin{center}
\begin{tabular}{llrlrr} \hline
\multicolumn{1}{c}{Identifier} &
\multicolumn{1}{c}{$resfinal$} &
\multicolumn{2}{c}{$totits$} &
\multicolumn{1}{c}{$maxbasis$} \\
\hline\Tstrut
Boeing/msc01050          &          7.68$\times 10^{-16}$  &       439  & (333)   & 162   \\ 
HB/bcsstk11              &          1.28$\times 10^{-14}$  &       602  & (644)   & 510   \\ 
HB/bcsstk26              &          3.21$\times 10^{-14}$  &       273  & (422)   & 148   \\ 
HB/bcsstk24              &          1.00$\times 10^{-14}$  &       755  & (834)   & 695   \\ 
HB/bcsstk16              &          2.05$\times 10^{-15}$  &        61 &  (80)   &  29   \\ 
Cylshell/s2rmt3m1        &          1.66$\times 10^{-15}$  &       262  & (473)   & 223   \\ 
Cylshell/s3rmt3m1        &          2.71$\times 10^{-15}$  &       841  & ($>$1000)   & 739   \\ 
Boeing/bcsstk38          &          4.01$\times 10^{-13}$  &       $>$2000    & ($>$1000) &    $>$1000   \\ 
Boeing/msc10848          &          4.36$\times 10^{-16}$  &       491 &  (622)   & 474   \\ 
Oberwolfach/t2dah$\_$e   &          2.35$\times 10^{-16}$  &        54 &  (34)   &  11   \\ 
Boeing/ct20stif          &          1.71$\times 10^{-12}$  &        $>$2000 &  ($>$1000)    & $>$1000  \\ 
DNVS/shipsec8            &          1.52$\times 10^{-12}$  &        $>$2000 &  ($>$1000)    & $>$1000  \\ 
Um/2cubes$\_$sphere      &          1.63$\times 10^{-16}$  &          26    &  (11) &     4     \\ 
GHS$\_$psdef/hood        &          2.13$\times 10^{-14}$  &           365  &  (346) &     189   \\ 
Um/offshore              &          3.32$\times 10^{-13}$  &        $>$2000 &  (101)   &   92  \\ 
\hline
\end{tabular}
\end{center}
}
\end{table}

\subsection{Results for $IC(\ell)$}
Figure~\ref{fig:changelevels} illustrates the influence of  the number of levels $\ell$
in the $IC(\ell)$ preconditioner computed in half precision and double precision.
Typically, $\ell$ is chosen to be small (large values lead to
slow computation times and loss of sparsity in $L$) although there are cases
where larger $\ell$ may be employed \cite{hypo:98}. In general, it is hoped that
as $\ell$ increases, the additional fill-in in $L$
will result in a better preconditioner (but there is no guarantee of this).
We observe that for the  well-conditioned problem UTEP/Dubcova2,
the half precision and double precision preconditioners generally behave 
in a similar way
(with the fp16 preconditioner having a higher iteration count for  $\ell=4$ and 5).
As $\ell$ increases from 1 to 9 the number of entries
in the incomplete factor increases from $1.76\times 10^6$ to $3.65\times 10^7$.
For the ill-conditioned problem Cylshell/s2rmt3m1,
the corresponding increase is from $1.43\times 10^5$ to $5.58\times 10^5$.
In this case, the iteration count for fp64+$IC(\ell)$
steadily decreases as $\ell$ increases but for fp16+$IC(\ell)$
the decrease is much less and it stagnates for $\ell > 4$.
In the rest of this section, we use $\ell = 3$.

\begin{figure}[htbp]
\begin{center}
\includegraphics[height=5cm]{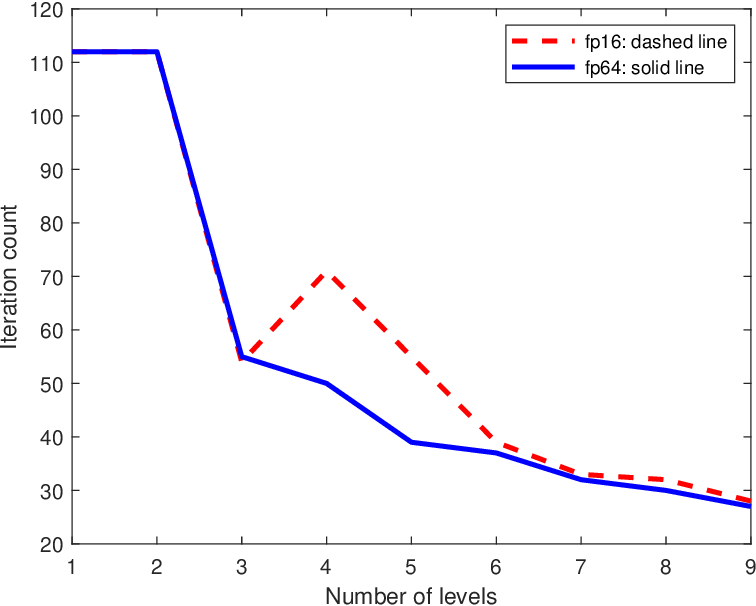} \hskip1.3cm
\includegraphics[height=5cm]{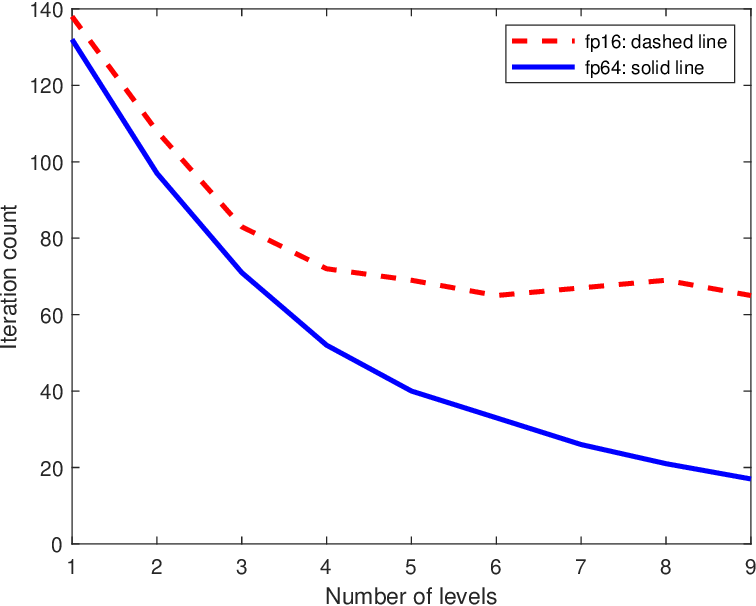} \\
\end{center}
\caption{IC-CG-IR total iteration counts
for the fp16+$IC(\ell)$ preconditioner (dashed line) and fp64+$IC(\ell)$ preconditioner
(solid line) with $l$ ranging from 1 to 9
for problems UTEP/Dubcova2 (left) and  Cylshell/s2rmt3m1 (right).} 
\label{fig:changelevels}
\end{figure}

\begin{table}[htbp]
\caption{Results for IC-CG-IR using an $IC(3)$ preconditioner: well-conditioned problems.
$resint$ and $resfinal$ are the initial and final scaled residuals.
$nnz(L)$ is the number of entries in the $IC(3)$ factor. $iouter$ and 
$totits$ denote the number of outer iterations and the total number of CG iterations, respectively.
$>1000$ indicates CG tolerance not reached on outer iteration $iouter$.
$nmod$ denotes the number of times problem B1 occurs during the factorization
(for fp64-$IC(3)$ it is equal to 0 for all our test cases and is omitted).
A count in bold indicates the fp16 result is within 10 per cent of (or is better than) the corresponding fp64 result.
}
\label{T:IC(3)+CGw}\vspace{3mm}
{
\footnotesize
\begin{center}
\begin{tabular}{lrlrrrr} \hline
\multicolumn{7}{c}{Preconditioner fp16-$IC(3)$} \\
{Identifier} &
{$resinit$} &
{$resfinal$} &
{$nnz(L)$} &
 {$iouter$} &
 {$totits$} &
 {$nmod$} \\
\hline\Tstrut
HB/bcsstk27              &     8.27$\times 10^{-5}$  &      1.23$\times 10^{-15}$  &      4.88$\times 10^4$ &     3  &    11   &     0 \\ 
Nasa/nasa2146            &     9.07$\times 10^{-5}$  &      7.41$\times 10^{-15}$  &      7.89$\times 10^4$ &     3  &    {\bf 11}    &     0 \\ 
Cylshell/s1rmq4m1        &     5.81$\times 10^{-5}$  &      1.05$\times 10^{-14}$  &      3.15$\times 10^5$ &     3  &    { 47}    &     0 \\ 
MathWorks/Kuu            &     1.74$\times 10^{-4}$  &      8.21$\times 10^{-15}$  &      7.65$\times 10^5$ &     3  &    {\bf 33}    &     0 \\ 
Pothen/bodyy6            &     7.39$\times 10^{-3}$  &      1.40$\times 10^{-16}$  &      1.76$\times 10^5$ &     4  &    108    &     2 \\ 
GHS$\_$psdef/wathen120   &     4.34$\times 10^{-4}$  &      6.69$\times 10^{-16}$  &      8.30$\times 10^5$ &     3  &     {\bf 6}    &     0 \\ 
GHS$\_$psdef/jnlbrng1    &     1.41$\times 10^{-3}$  &      5.10$\times 10^{-15}$  &      2.77$\times 10^5$ &     3  &    {\bf 12}    &     0 \\ 
Williams/cant            &     3.28$\times 10^{-4}$  &      3.74$\times 10^{-15}$  &      9.95$\times 10^6$ &     3  &  { 1103}    &     0 \\ 
UTEP/Dubcova2            &     1.54$\times 10^{-3}$  &      1.88$\times 10^{-13}$  &      6.22$\times 10^6$ &     3  &    {\bf 54}    &     0 \\ 
Cunningham/qa8fm         &     3.54$\times 10^{-4}$  &      1.04$\times 10^{-16}$  &      5.14$\times 10^6$ &     3  &     {\bf 6}    &     0 \\ 
Mulvey/finan512          &     3.50$\times 10^{-4}$  &      4.97$\times 10^{-17}$  &      4.08$\times 10^6$ &     3  &     {\bf 6}    &     0 \\ 
GHS$\_$psdef/apache1     &     8.51$\times 10^{-5}$  &      3.65$\times 10^{-14}$  &      1.54$\times 10^6$ &     2  &   {\bf 116}    &     0 \\ 
Williams/consph          &     1.04$\times 10^{-4}$  &      1.57$\times 10^{-14}$  &      2.02$\times 10^7$ &     3  &   {\bf 160}    &     0 \\ 
AMD/G2$\_$circuit        &     8.15$\times 10^{-4}$  &      3.28$\times 10^{-16}$  &      1.04$\times 10^6$ &     4  &   {\bf 282}    &     0 \\ 
\hline\hline
\multicolumn{7}{c}{Preconditioner fp64-$IC(3)$} \\
{Identifier} &
{$resinit$} &
{$resfinal$} &
{$nnz(L)$} &
 {$iouter$} &
 {$totits$} &
 {}  \\
\cline{1-6}\Tstrut
HB/bcsstk27              &     3.93$\times 10^{-5}$  &      3.13$\times 10^{-16}$  &      4.88$\times 10^4$ &     3  &     8    &      \\ 
Nasa/nasa2146            &     4.32$\times 10^{-5}$  &      8.80$\times 10^{-16}$  &      7.89$\times 10^4$ &     3  &    11    &      \\ 
Cylshell/s1rmq4m1        &     3.55$\times 10^{-5}$  &      1.66$\times 10^{-15}$  &      3.15$\times 10^5$ &     3  &    56    &      \\ 
MathWorks/Kuu            &     2.50$\times 10^{-4}$  &      2.25$\times 10^{-15}$  &      7.65$\times 10^5$ &     3  &    30    &      \\ 
Pothen/bodyy6            &     9.31$\times 10^{-3}$  &      1.40$\times 10^{-16}$  &      1.76$\times 10^5$ &     4  &    63    &      \\ 
GHS$\_$psdef/wathen120   &     1.77$\times 10^{-4}$  &      9.56$\times 10^{-17}$  &      8.30$\times 10^5$ &     3  &     6    &        \\ 
GHS$\_$psdef/jnlbrng1    &     9.40$\times 10^{-4}$  &      3.22$\times 10^{-16}$  &      2.77$\times 10^5$ &     3  &    12    &         \\ 
Williams/cant            &     3.19$\times 10^{-4}$  &      6.09$\times 10^{-15}$  &      9.95$\times 10^6$ &     3  &   927    &      \\ 
UTEP/Dubcova2            &     1.65$\times 10^{-3}$  &      1.84$\times 10^{-13}$  &      6.22$\times 10^6$ &     3  &    55    &      \\ 
Cunningham/qa8fm         &     3.29$\times 10^{-5}$  &      1.24$\times 10^{-16}$  &      5.14$\times 10^6$ &     3  &     5    &      \\ 
Mulvey/finan512          &     3.72$\times 10^{-5}$  &      6.63$\times 10^{-17}$  &      4.08$\times 10^6$ &     3  &     6    &      \\ 
GHS$\_$psdef/apache1     &     8.28$\times 10^{-5}$  &      4.12$\times 10^{-14}$  &      1.54$\times 10^6$ &     2  &   120    &         \\ 
Williams/consph          &     3.39$\times 10^{-4}$  &      6.07$\times 10^{-3}$   &      2.02$\times 10^7$ &     1  &  $>1000$    &      \\ 
AMD/G2$\_$circuit        &     8.42$\times 10^{-5}$  &      2.46$\times 10^{-16}$  &      1.04$\times 10^6$ &     4  &   273    &        \\ 
\cline{1-6}
\end{tabular}
\end{center}
}
\end{table}

\begin{table}[htbp]
\caption{Results for IC-CG-IR and IC-GMRES-IR using an $IC(3)$ preconditioner: ill-conditioned problems.
$resint$ is the initial scaled residual;
$resfinal$ is the final IC-CG-IR scaled residual. $nnz(L)$ is the number of entries in the $IC(3)$ factor. $iouter$ and 
$totits$ denote the number of outer iterations and the total number of CG iterations with the GMRES
statistics in parentheses.
$>1000$ indicates CG (or GMRES) tolerance not reached on outer iteration $iouter$.
$nmod$ and $nofl$ denote the numbers of times problems B1 and $B3$ occur during the factorization
(for fp64-$IC(3)$ they are equal to 0 for all our test cases and is omitted).
A count in bold indicates the fp16 result is within 10 per cent (or better) of the corresponding fp64 result.
\ddag~indicates failure to compute 
the factorization because of enormous growth in its entries.
}
\label{T:IC(3)+CGi}\vspace{3mm}
{
\footnotesize
\begin{center}
\begin{tabular}{llllrrlrr} \hline
\multicolumn{8}{c}{Preconditioner fp16-$IC(3)$} \\
\multicolumn{1}{c}{Identifier} &
\multicolumn{1}{c}{$resinit$} &
\multicolumn{1}{c}{$resfinal$} &
\multicolumn{1}{c}{$nnz(L)$} &
\multicolumn{1}{c}{$iouter$} &
\multicolumn{2}{c}{$totits$} &
\multicolumn{1}{c}{$nmod$} &
\multicolumn{1}{c}{$nofl$} \\
\hline\Tstrut
Boeing/msc01050          &     1.61$\times 10^{-3}$  &      6.47$\times 10^{-14}$  &      3.74$\times 10^4$ &     3  &   125 & (64)   &     2 &     0 \\ 
HB/bcsstk11              &     5.84$\times 10^{-4}$  &      1.89$\times 10^{-13}$  &      4.18$\times 10^4$ &     3  &   265 & (184)   &     0 &     2 \\ 
HB/bcsstk26              &     6.42$\times 10^{-5}$  &      1.51$\times 10^{-13}$  &      3.43$\times 10^4$ &     3  &    80 & (70)   &     2 &     0 \\ 
HB/bcsstk24              &     5.91$\times 10^{-7}$  &      1.77$\times 10^{-13}$  &      2.27$\times 10^5$ &     3  &   437 & (260)   &     0 &     2 \\ 
HB/bcsstk16              &     6.16$\times 10^{-4}$  &      6.64$\times 10^{-15}$  &      4.89$\times 10^5$ &     3  &    17 &  (17)  &     0 &     0 \\ 
Cylshell/s2rmt3m1        &     9.24$\times 10^{-6}$  &      4.11$\times 10^{-15}$  &      2.60$\times 10^5$ &     3  &    83 & (117)   &     0 &     0 \\ 
Cylshell/s3rmt3m1        &     1.72$\times 10^{-6}$  &      9.01$\times 10^{-15}$  &      2.60$\times 10^5$ &     3  &   {\bf 386} &  (504)  &     1 &  1 \\ 
Boeing/bcsstk38          &     1.07$\times 10^{-3}$  &      1.43$\times 10^{-15}$  &      5.64$\times 10^5$ &     4  &   $1004$ &  (282)  &     2 &     0 \\ 
Boeing/msc10848          &     5.97$\times 10^{-7}$  &      1.13$\times 10^{-14}$  &      2.51$\times 10^6$ &     3  &   138 & (89)   &     0 &     0 \\ 
Oberwolfach/t2dah$\_$e   &     6.81$\times 10^{-4}$  &      1.88$\times 10^{-16}$  &      3.29$\times 10^5$ &     3  &     {\bf 6} & ({\bf 6})   &     0 &     0 \\ 
Boeing/ct20stif          &     2.78$\times 10^{-5}$  &      1.63$\times 10^{-9}$   &      6.70$\times 10^6$ &     3  &  $>1000$ & ($>1000$)   &     2 &     0 \\ 
DNVS/shipsec8            &     5.56$\times 10^{-6}$  &      1.26$\times 10^{-16}$  &      1.22$\times 10^7$ &     4  &  $2067$ &  ($1399$)   &     2 &     0 \\ 
Um/2cubes$\_$sphere      &     1.02$\times 10^{-3}$  &      1.63$\times 10^{-16}$  &      8.70$\times 10^6$ &     3  &     {\bf 6} & ({\bf 6})   &     0 &     0 \\ 
GHS$\_$psdef/hood        &     6.80$\times 10^{-4}$  &      5.01$\times 10^{-17}$  &      2.78$\times 10^7$ &     4  &   {\bf 444} &  ({\bf 406})  &     0 &     4 \\ 
Um/offshore              &     1.53$\times 10^{-3}$  &      1.38$\times 10^{-13}$  &      2.08$\times 10^7$ &     3  &   {\bf 103} & ({\bf 40})   &     0 &     4 \\ 
\hline
\multicolumn{8}{c}{Preconditioner fp64-$IC(3)$} \\
\multicolumn{1}{c}{Identifier} &
\multicolumn{1}{c}{$resinit$} &
\multicolumn{1}{c}{$resfinal$} &
\multicolumn{1}{c}{$nnz(L)$} &
\multicolumn{1}{c}{$iouter$} &
\multicolumn{2}{c}{$totits$} &
\multicolumn{1}{c}{} &
\multicolumn{1}{c}{ } \\
\cline{1-7}\Tstrut
Boeing/msc01050          &     1.41$\times 10^{-4}$  &      5.80$\times 10^{-14}$  &      3.74$\times 10^4$ &     3  &    38 & (25)    &     &      \\ 
HB/bcsstk11              &     1.48$\times 10^{-5}$  &      6.96$\times 10^{-14}$  &      4.18$\times 10^4$ &     3  &    29 & (29)   &     &      \\ 
HB/bcsstk26              &     9.24$\times 10^{-5}$  &      1.68$\times 10^{-13}$  &      3.43$\times 10^4$ &     3  &    60 & (62)   &     &      \\ 
HB/bcsstk24              &     4.01$\times 10^{-7}$  &      1.01$\times 10^{-13}$  &      2.27$\times 10^5$ &     3  &    71 &  (67)  &     &      \\ 
HB/bcsstk16              &     7.01$\times 10^{-4}$  &      3.13$\times 10^{-15}$  &      4.89$\times 10^5$ &     3  &    15 &  (14)  &     &      \\ 
Cylshell/s2rmt3m1        &     8.61$\times 10^{-6}$  &      8.74$\times 10^{-15}$  &      2.60$\times 10^5$ &     3  &    71 &  (105)  &     &      \\ 
Cylshell/s3rmt3m1        &     2.00$\times 10^{-6}$  &      6.12$\times 10^{-5}$   &      2.60$\times 10^5$ &     1  &  $>1000$ &  ($>1000$)  &     &      \\ 
Boeing/bcsstk38          &     4.67$\times 10^{-9}$  &      8.34$\times 10^{-14}$  &      5.64$\times 10^5$ &     3  &   154 &  115)  &     &      \\ 
Boeing/msc10848          &     8.73$\times 10^{-10}$ &      2.26$\times 10^{-16}$  &      2.51$\times 10^6$ &     3  &    47 &  (46)  &     &      \\ 
Oberwolfach/t2dah$\_$e   &     5.49$\times 10^{-6}$  &      6.54$\times 10^{-16}$  &      3.29$\times 10^5$ &     3  &     5 &  (5)  &     &   \\ 
Boeing/ct20stif          &     5.82$\times 10^{-8}$  &      2.02$\times 10^{-11}$  &      6.70$\times 10^6$ &     3  &  $>1000$ & ($>1000$)  &     &      \\ 
DNVS/shipsec8            &     7.50$\times 10^{-7}$  &      1.11$\times 10^{-16}$  &      1.22$\times 10^7$ &     4  &   313 &  (181)  &     &      \\ 
Um/2cubes$\_$sphere      &     2.25$\times 10^{-5}$  &      1.63$\times 10^{-16}$  &      8.70$\times 10^6$ &     3  &     5  &  (5) &     &   \\ 
GHS$\_$psdef/hood        & \ddag & \ddag   &      \ddag  &     \ddag  &   \ddag    & \ddag  &   &      \\ 
Um/offshore              & \ddag &  \ddag  &      \ddag  &     \ddag  &   \ddag    & \ddag  &    &      \\ 
\cline{1-7}
\end{tabular}
\end{center}
}
\end{table}

Tables~\ref{T:IC(3)+CGw} and \ref{T:IC(3)+CGi} present results for IC-CG-IR with the fp16-$IC(3)$ and fp64-$IC(3)$ preconditioners
for the well-conditioned and ill-conditioned test sets, respectively. The latter also reports total iteration counts for IC-GMRES-IR
(the number of outer iterations $iouter$ for IC-GMRES-IR and IC-CG-IR are the same for all the test examples).
B3 breakdowns occur for a small number of the ill-conditioned problems.
We see that, for our examples, the number of entries in $L$ is (approximately)
the same for both fp16 and fp64 arithmetic, indicating that (in contrast
to $IC(0)$) the number of subnormal numbers is small.
fp16-$IC(3)$ performs as well as fp64-$IC(3)$ on most of the
well-conditioned problems. Note that for problem Williams/consph
 the required accuracy was not achieved using fp64 on the first outer iteration. 
For the ill-conditioned problems, while the fp16-$IC(3)$ preconditioner
combined with IC-CG-IR and IC-GMRES-IR is able to return a computed solution with a small residual, for the problems
for which B1 and/or B3 are detected in fp16 arithmetic the 
iteration counts are significantly greater than for the fp64-$IC(3)$ preconditioner (see, for instance,
Boeing/msc01050 and HB/bcsstk11). This is because the use of shifts to prevent breakdowns means the 
computed factors are for a shifted matrix. 
However, for a small number of problems the computation of the fp64-$IC(3)$ preconditioner
fails because the computed factor has very large entries, which do not overflow
in double precision but make it useless as a preconditioner.
This indicates that it is not just in half precision arithmetic that it is necessary to monitor the 
possibility of growth occurring in the factor entries, something that is not currently considered when
computing incomplete Cholesky factorizations in double (or single) precision arithmetic.

Note that, because the $L$ factor computed using
fp16 arithmetic is less accurate than that computed using fp64 arithmetic,
the initial residual scaled $resint$ 
given by (\ref{eq:resid}) with $x = L^{-T}L^{-1}$ is typically
larger for fp16  than for fp64 arithmetic. However, if double precision
accuracy in the computed solution is not required,
performing refinement may be unnecessary
(or a small number of steps may be sufficient),
even for fp16 incomplete factors. 

\section{Concluding remarks and future directions}\label{sec:conclusions}

The focus of this study is the construction 
and employment of incomplete factorizations using half precision arithmetic
to solve large-scale sparse linear systems to double precision accuracy.
Our experiments, which simulate half precision arithmetic through the use of the NAG compiler
demonstrate that, when carefully implemented, the use of fp16  
level-based incomplete factorization preconditioners may not
impact on the overall accuracy of the computed solution, even when
a small tolerance is imposed on the requested scaled residual. Unsurprisingly, the number of iterations of
the Krylov subspace method that is used in the refinement process can be
greater for fp16 factors compared to fp64 factors
but generally this increase is only significant for highly ill-conditioned systems.
Our results support the view 
that, for many real-world problems, it is sufficient to employ half precision arithmetic.
Its use may be particularly advantageous if the linear system does not need to
be solved to high accuracy. Our study also encourages us to conjecture that by building safe operations
into sparse direct solvers it should be possible
to build efficient and robust half precision variants and,
because this would lead to substantial memory savings for the matrix factors, it
could potentially  allow direct solvers to be used (in combination with 
an appropriate refinement process) to solve much larger problems than is currently possible. {This is a future direction
that is becoming more feasible to
explore as compiler support for fp16 arithmetic
improves and becomes available on more platforms
and for more languages.}

It is of interest to consider other Krylov solvers
such as flexible CG and flexible GMRES  and to explore other classes of algebraic preconditioners,
to consider how they can be safely computed using fp16 arithmetic and how effective
they are compared to higher precision versions.
For sparse approximate inverse (SPAI) preconditioners,  
we anticipate that it may be possible to combine the systematic dropping of subnormal
quantities with avoiding overflows in local linear solves as we
have discussed without 
significantly affecting the preconditioner quality because such changes may only influence
a small number of the computed columns of the SPAI preconditioner
(see also \cite{cakh:2023} for a recent analysis of SPAI preconditioners in mixed precision). 
For other approximate inverses, such as AINV \cite{bemt:96} 
and AIB \cite{saad:03}, the sizes of the diagonal entries follow from 
maintaining a generalized orthogonalization property and avoiding overflows 
using local or global modifications may be possible but challenging.
We also plan to consider the construction of the 
{\tt HSL\_MI28} \cite{sctu:2014a} IC preconditioner using low precision arithmetic.
{\tt HSL\_MI28} uses an extended memory approach for the 
robust construction of the factors.
Once the factors have been computed, the extended memory can be freed,
limiting the size of the factors and the work needed in the substitution steps, generally
without a significant reduction in the preconditioner quality.
A challenge here is safely allowing
intermediate quantities of absolute value at least $x_{min}$ to be retained in the factors or 
in the extended memory.

\medskip
\noindent
{\bf Acknowledgements} We are very grateful to the anonymous reviewers for their detailed
comments that have led to significant improvements to this paper.

%\bibliographystyle{ACM-Reference-Format}
%\bibliography{btbook}
\def\cprime{$'$} \def\cprime{$'$} \def\cprime{$'$}

\end{document}